\documentclass[12pt]{article}
\usepackage{acronym}
\usepackage{amsfonts}
\usepackage{amsmath,amssymb,bm}
\usepackage{graphicx}
\usepackage{mdframed}
\usepackage{multirow}
\usepackage{epstopdf}
\usepackage{tabularx}
\usepackage{float}

\newtheorem{theorem}{Theorem} 

\newtheorem{lemma}{Lemma}[section]
\newtheorem{definition}{Definition}[section]

\def\be{\begin{Example}}
\def\ee{\end{Example}}
\def\bt{\begin{theorem}}
\def\et{\end{theorem}\bigskip}
\def\bl{\begin{Lemma}}
\def\el{\end{Lemma}\bigskip}
\def\ep{\end{Proposition}\bigskip}
\def\bp{\begin{Proposition}}
\def\bd{\begin{definition}}
\def\ed{\end{definition}}

\begin{document}
\title{\bf Accelerating Power Methods for Higher-order Markov Chains
}
\author{
Gaohang Yu\footnotemark[1], Yi Zhou\footnotemark[1]~, Laishui Lv\footnotemark[2]~
}
\renewcommand{\thefootnote}{\fnsymbol{footnote}}

\footnotetext[1]{ School of Science, Hangzhou dianzi University, 310018, China. E-mail: maghyu@163.com(G.Yu) \ 13592051764@163.com(Y.Zhou)}
\footnotetext[2]{School of Computer Science and Engineering, Nanjing University of Science and Technology, Nanjing, 210094, China. E-mail: malslv@163.com}
\maketitle
\begin{abstract}
Higher-order Markov chains play a very important role in many fields, ranging from multilinear PageRank to financial modeling. In this paper, we propose three accelerated higher-order power methods for computing the limiting probability distribution of higher-order Markov chains, namely higher-order power method with momentum and higher-order quadratic extrapolation method. The convergence results are established, and numerical experiments are reported to show the efficiency of the proposed algorithms. In particular, the non-parametric quadratic extrapolation method is very competitive, and outperforms state-of-the-art competitions.

\end{abstract}
\section{Introduction}

Markov chains are powerful tools to analyze and predict traffic flows, communications networks, genetic issues, and a variety of stochastic (probabilistic) processes over time, in which the probability of each event depends only on the state attained in the previous event. Considering a stochastic process $\{X_t, t = 0, 1, 2, \ldots \}$ that
takes on a finite set $\{1, 2,\ldots, n\}\equiv \ \langle n \rangle$. An element in $\langle n \rangle$ is called a state of the process.
The definition of a Markov Chain can be given as follows.

\begin{definition}
Assume there exists a fixed probability $p_{i,j}$ independent of time such that
$$Prob(X_{t+1}=i|X_t=j,X_{t-1}=i_{t-1},\ldots,X_0=i_0)=Prob(X_{t+1}=i|X_t=j)=p_{i,j},$$
where $i,j,i_0,i_1,\ldots,i_{t-1}\in \langle n \rangle$ and $\{X_t\}(t=0,1,2,\ldots)$ is a stochastic process. Then this is called a Markov chain process.
\end{definition}

The probability $p_{i,j}$ represents the probability that the process will make a transition
to state $i$ given that currently the process is state $j$. Clearly one has
$$p_{i,j}\ge0,\ \sum^n_{i=1}p_{i,j} = 1, j = 1,\ldots, n.$$
The matrix $P = (p_{i,j})$ is called the one-step transition probability matrix of the
process. A vector $\bar{x}$ is said to be a limiting or stationary probability distribution of a
finite Markov chain having $n$ states with
$$\bar{x}_i \ge 0,\ \forall i, \ \sum^n_{i=1}\bar{x}_i = 1, \mbox{and} P\bar{x} = \bar{x}.$$

In real world, there are many situations that one would like to employ higher-order Markov
chain models as a mathematical tool to analyze data sequences, in which the probability of $X_{t+1}=i$ not only depends on the adjacent time state $X_t$ but also depends on more previous time states. The $(m-1)^{th}$ order Markov chain model is given as follows.

\begin{definition}
Assume there exists a fixed probability $p_{i_1,i_2,\ldots,i_m}$ independent of time such that
$$0\leq p_{i_1,i_2,\ldots,i_m}=Prob(X_{t+1}=i_1|X_t=i_2,\ldots,X_{t-m+2}=i_m)\leq1,$$
where $i_1,\ldots,i_{m}\in \langle n \rangle$ and $\sum_{i_1=1}^np_{i_1,i_2,\ldots,i_m}=1$. Then this is called a $(m-1)^{th}$ order Markov chain process.
\end{definition}

It is clear that the $(m-1)^{th}$ order Markov chain process will reduces to first-order Markov chain when $m=2$. The probability $p_{i_1,i_2,\ldots,i_m}$ represents that process make transition to the state $i_1$ given that currently the process is in the state $i_2$ and previously the process is in the states $i_3,\ldots,i_m$. Tensor $\mathcal{P}=(p_{i_1,i_2,\ldots,i_m})$ is called transition probability tensor. A number of applications can be found in the literature, for example, chemistry\cite{FSJ84,KKS11}, physics\cite{ANM07} and multilinear PageRank\cite{Gleich14}.

In \cite{WMNG14}, Li and Ng established the following approximated tensor model for the stationary probability distribution of Higher-order Markov chains:
\begin{equation}\label{approximated model}
 x=\mathcal{P}x^{m-1},\ \  x\geq0,\ \ \|x\|_1=1,
\end{equation}
where $\mathcal{P}x^{m-1}$ is defined by:
$$
(\mathcal{P}x^{m-1})_i=\sum_{i_2,\ldots,i_m=1}^n p_{i i_2 \cdots i_m}x_{i_2}\cdots
x_{i_m},\ \ i=1,2,\ldots,n,
$$
and $x=(x_i)$ is called a stationary probability distribution vector of higher-order Markov chains. 
The stationary probability distribution vector is unique under suitable conditions \cite{Chang13,Geiger17,WMNG14,HuangQi18}.

Later, many researchers employed the higher-order Markov chains model to explore some applications such as in random walk\cite{Gleich16} and
multilinear PageRank\cite{Gleich14}. Gleich, Lim and Yu \cite{Gleich14} first studied the following multilinear PageRank model:
\begin{equation}\label{multiPageRank}
 x=\theta\hat{\mathcal{P}}x^{m-1}+(1-\theta)v,
\end{equation}
where tensor $\hat{\mathcal{P}}$ is a transition probability tensor, $v$ is transition probability vector, and $\theta\in(0,1)$ is a damping parameter. We can rewrite the equation  (\ref{multiPageRank}) as follows
\begin{equation}\label{remultiPageRank}
 x=\mathcal{P}x^{m-1},\ \|x\|_1=1, \ \mathcal{P}=\theta\hat{\mathcal{P}}+(1-\theta)*\mathcal{V},
\end{equation}
where $\mathcal{V}=(v_{i_1i_2\ldots,i_m})$ with $v_{i_1i_2\ldots,i_m}=v_{i_1},\forall i_2,\ldots,i_m$. It is easy to see that the tensor $\mathcal{P}$ is also a transition probability tensor.

Recently, Li {\it et al.} \cite{LiLiu20} investigated the uniqueness of the fixed-point for the equation (\ref{remultiPageRank}) and presented perturbation analysis for Multilinear PageRank model (\ref{multiPageRank}). In \cite{Gleich14}, several iterative algorithms (a fixed-point algorithm, a shifted fixed-point algorithm, a inner-outer iteration algorithm, an inverse iteration algorithm and a Newton algorithm) are proposed by Gleich {\it et al.}. Furthermore, Meini and Poloni \cite{Meini18} proposed the Perron-based iteration and Cipolla {\it et al.} \cite{Tudisco20} presented some extrapolation methods for fixed-point multilinear PageRank computations.

 As for solving the tensor equations (\ref{approximated model}), Li and Ng\cite{WMNG14} extended the power method to compute the tensor equation (\ref{approximated model}). They given the convergence analysis of the proposed iterative algorithm. In \cite{Liu19}, Liu {\it et al.} proposed several relaxation algorithms for solving equation (\ref{approximated model}). And a truncated power method is presented in \cite{DingNgWei18} for sparse Markov chains.
Power-type methods are very popular due to their simplicity and efficiency, especially for large-scale problems \cite{DingNgWei18,Gleich14}. However, as shown in \cite{WuChu17}, the convergence rate of the higher-order power method will be slow when the spectral gap is small. Moreover, they point out that there exists irreducible and aperiodic transition probability tensors where the Z-eigenvector type power iteration fails to converge.

A pair $(\lambda,x)\in \mathbb{R}\times\mathbb{R}^n\backslash\{0\}$ is called $Z$-eigenpair of the tensor $\mathcal{P}$ if
\begin{equation}\label{Zeigenpair}
\mathcal{P}x^{m-1}=\lambda x,\ \ \|x\|=1.
\end{equation}
This definition was proposed by Qi\cite{Qi05} and Lim\cite{Lim05}, independently. Here, $\|\cdot\|$ could be $l_1$-norm or $l_2$-norm. Moreover, for $x\in R^n$ with $\|x\|_1=1$, $(x,\lambda)$ is the $Z_1$-eigenpair if and only if $(\frac{x}{\|x\|_2},\frac{\lambda}{\|x\|_2^{m-1}})$ is the corresponding $Z_2$-eigenpair \cite{Chang13}. It is clear that the solution to the tensor equations (\ref{approximated model}), i.e. the stationary probability distribution vector $x$, is the $Z_1$-eigenvector of $\mathcal{P}$ while 1 is the largest $Z_1$-eigenvalue.

In this paper, we propose three algorithms for solving the stationary probability distribution of higher-order Markov chains by accelerating the higher-order power method. The main contributions of this paper are
\begin{itemize}
  \item to present two novel higher-order power methods with momentum for Z-eigenvector computations of tensor;

  \item to propose a non-parametric higher-order quadratic extrapolation method to compute the stationary probability distribution of higher-order Markov chains;

  \item to establish the convergence theorems for the proposed algorithms;

  \item to use the proposed algorithms for some applications such as fixed-point multilinear PageRank computations.
\end{itemize}

\par The rest of this paper is organized as follows. Firstly, some preliminary knowledge and existing methods are presented in Section 2. In Section 3, we propose three provable power-type methods for calculating the limiting probability distribution vector of higher-order Markov chains. In Section 4, the convergence theorems for the proposed methods are established. Numerical experiments are given and analyzed in Section 5. The last section is the conclusions.

Before end of this section, we would like to describe notations and show some preliminary knowledge on tensors.For the details of basic definitions and properties of tensors, we refer to reader to\cite{QiLuo17} and the references therein.

Let $\mathbb{R}$ be the real field. An $m$th-order $n$-dimensional real tensor $\mathcal{P}$ consists of $n^m$ entries in real numbers:
$
\mathcal{P}=(a_{i_1i_2\cdots
i_m}),\,\, p_{i_1i_2\cdots i_m}\in \mathbb{R},\,\, \mbox{for any}  \ \ i_1,i_2,\ldots,i_m \in [n],
$
where $[n]=\{1,2,\ldots,n\}$. $\mathcal{P}$ is called non-negative(or,respectively,positive) if $p_{i_1i_2\cdots
i_m}\geq0$(or, respectively, $p_{i_1i_2\cdots i_m}>0$).
Given two vectors $x,y\in R^n$, we define
$$
\mathcal{P}(x^{m-1}-y^{m-1})\equiv \mathcal{P}x^{m-1}-\mathcal{P}y^{m-1}.
$$

Let $x_+=max(x,0)$ and $proj(x)=\frac{x_+}{\|x_+\|_1}$. It is easy to get that $proj(x)$ is a transition probability vector. We also show the definition of irreducible tensors as follows.
\begin{definition}\label{dy1}
An $m$-order $n$-dimensional tensor $\mathcal{P}$ is called reducible if there exists
a nonempty proper index subset $I\subset\{1,2,\ldots,n\}$ such that
$$P_{i_1i_2\ldots i_m}=0,\ \ \forall i_1\in I,\ \ \forall i_2,\ldots,i_m\notin I.$$
If $\mathcal{P}$ is not reducible, then we call $\mathcal{P}$ irreducible.
\end{definition}

\section{The existing methods}

 To some extent, the limit probability distribution problem (\ref{approximated model}) can be transformed into an optimization problem:
\begin{equation}\label{min-exp}
\min f(x)=\min_{\|x\|_1=1}\frac{1}{2}x^Tx-\frac{1}{m}x^T \mathcal{P}x^{m-2}x,
\end{equation}
where $\mathcal{P}x^{m-2}$ is a matrix with its component as $$(\mathcal{P}x^{m-2})_{i_1i_2}=\sum_{i_3,\ldots,i_m=1}^n p_{i_1 i_2i_3 \cdots i_m}x_{i_3}\cdots
x_{i_m}\ \ \mbox{for all}\ \ i_1,i_2=1,\cdots,n.$$
Let $\nabla f(x)=x-\mathcal{P}x^{m-1}$, and minimizing the above problems via gradient descent scheme

\begin{equation*}
  \begin{split}
     x_k&=x_{k-1}-\alpha \nabla f(x_{k-1}) \\
        &=x_{k-1}-\alpha(x_{k-1}-\mathcal{P}x_{k-1}^{m-1}),
   \end{split}
\end{equation*}
setting $\alpha=1$, then we can get the following higher-order power method(HOPM) for solving the tensor equation (\ref{approximated model}).

\noindent\rule{\textwidth}{1pt}
\underline{\textbf{Algorithm 1:higher-order power method(HOPM) \cite{WMNG14}}} \\
1. Given a transition probability tensor $\mathcal{P}$, maximum $k_{max}$, termination tolerance $\epsilon$ and an initial point $x_0$;\\
2. Initialize $k=1$.\\
3.\ \ \ \ $x_{k}=\mathcal{P}x_{k-1}^{m-1}$;\\
4. \ \ \  $\delta=\|x_{k}-x_{k-1}\|$;\\
5. \ \ \  $k=k+1$;\\
6. \textbf{until} $\delta<\epsilon$.\\
\noindent\rule{\textwidth}{1pt}
\textbf{Remark 1}. The main computational cost of the algorithm depends on the cost of performing tensor operation. Assume that there are $O(N)$ nonzero entries (sparse
data) in tensor $\mathcal{P}$, the cost of this tensor calculation are of $O(N)$ arithmetic operations. Under some suitable conditions, they established
the linear convergence of the above algorithm.

In \cite{KoldaMayo14}, Kolda and Mayo presented an adaptive, monotonically convergent, shifted power method for computing  tensor Z-eigenpairs, called GEAP method.

\noindent\rule{\textwidth}{1pt}
\underline{\textbf{Algorithm 2: Z-Eigenpair Adapative Power Method (GEAP Method) \cite{KoldaMayo14} }} \\
Given a transition probability tensor $\mathcal{P}$, maximum $k_{max}$, termination tolerance $\epsilon$ and an initial point $x_0$; Let $\tau>0$ is the tolerance on being positive definite.\\
Initialize $k=1$.\\
\textbf{1.} Precompute $\mathcal{P}x_{k-1}^{m-2}$,$\mathcal{P}x_{k-1}^{m-1}$\\
\textbf{2.} $H_{k-1}\leftarrow H(x_{k-1})=m(m-1)\mathcal{P}x_{k-1}^{m-2}$\\
\textbf{3.} $\alpha_k\leftarrow \max\{0,(\tau-\lambda_{min}( H_{k-1}))/m\}$\\
\textbf{4.} $\hat{x}_{k}\leftarrow \mathcal{A}x_{k-1}^{m-1}+\alpha_kx_{k-1} $, $x_{k}=proj(\hat{x}_{k})$.\\
\textbf{5.} $\delta=\|x_{k}-x_{k-1}\|$;\\
\textbf{6.}   $k=k+1$;\\
\textbf{7.}\textbf{until} $\delta<\epsilon$.\\
\noindent\rule{\textwidth}{1pt}

In \cite{Liu19}, Liu et al. proposed several relaxation methods for computing tensor equation (\ref{approximated model}). In particular, by using relaxation technique to the higher-order power method, they developed a novel algorithm as follows.

\noindent\rule{\textwidth}{1pt}
\underline{\textbf{Algorithm 3: Relaxation higher-order power method, RHOPM (Alg.2 in \cite{Liu19})}} \\
1. Given a transition probability tensor $\mathcal{P}$, $\gamma>0$, termination tolerance $\epsilon$ and an initial point $x_0$;\\
2. Initialize $k=1$.\\
3.\ \ \ \ $y_{k}=\mathcal{P}x_{k-1}^{m-1}$;\\
4.\ \ \ \ $\hat{x}_{k}=\gamma y_k+(1-\gamma)x_{k-1}$, $x_k=proj(\hat{x}_{k})$\\
5. \ \ \  $\delta=\|x_{k}-x_{k-1}\|$;\\
6. \ \ \  $k=k+1$;\\
7. \textbf{until} $\delta<\epsilon$.\\
\noindent\rule{\textwidth}{1pt}


\section{Accelerating Power methods}
In this section, we will propose three accelerated higher-order power methods, including two higher-order power methods with momentum term and quadratic extrapolation higher-order power method, respectively.
\subsection{Higher-order power methods with momentum}
In \cite{Polyak04,Nesterov04}, some accelerated first-order methods are proposed by adding momentum terms to classic gradient method, called heavy-ball method and Nesterov's accelerated gradient method (NAG), respectively. Recently, Xu {\it et al.} \cite{xu2018} proposed a power method with momentum for principal component analysis.

Along this line, motivated by the efficiency of momentum methods, we propose the following two algorithms for solving the tensor equations (\ref{approximated model}) by respectively adding two different momentum terms to higher-order power method, referred to as the HOPMM-I and HOPMM-II.

\noindent\rule{\textwidth}{1pt}
\underline{\textbf{HOPMM-I}} \\
1. Given a transition probability tensor $\mathcal{P}$, maximum $k_{max}$, $\beta>0$, termination tolerance $\epsilon$ and an initial point $x_0$;\\
2. Initialize $k=1$.\\
3.  \textbf{repeat}\\
4.\ \ \ \ $x_{k}=\mathcal{P}x_{k-1}^{m-1}$;\\
5. periodically,\\
6. \ \ \ \ \ \ $\hat{x}_{k}=x_{k}+\beta(x_{k-1}-x_{k-2})$, $x_{k}=proj(\hat{x}_{k})$;\\
7: If $\|x_{k}-x_{k-1}\|_1<\epsilon$, break and output $x_k$;\\
8. $k=k+1$, back to step 4.\\
\noindent\rule{\textwidth}{1pt}
\textbf{Remark} 2. The $\beta(x_{k-1}-x_{k-2})$ is called momentum term. By choosing a suitable parameter $\beta$, the HOPMM-I will performs better than higher-order power method. In particular, if $\beta=0$, the HOPMM-I will reduces to higher-order power method that proposed by Li et al.\cite{WMNG14}. Compared with the RHOPM, our proposed method uses three iterative points to generate next iterative point. In the HOPMM-I, we will execute the momentum extrapolation at every 3 steps.

Furthermore, we can also add a ``heavy-ball" momentum term to higher-order power method and obtain the following algorithm.

\noindent\rule{\textwidth}{1pt}
\underline{\textbf{HOPMM-II}} \\
1: Given a transition probability tensor $\mathcal{P}$, maximum $k_{max}$, $\eta>0$, termination tolerance $\epsilon$ and an initial point $x_0$;\\
2:Initialize $k=1$.\\
3.  \textbf{repeat}\\
4: \ \ \ \ $x_k=\mathcal{P}x_{k-1}^{m-1}$;\\
5. periodically,\\
6: \ \ \ \ $\hat{x}_{k}=x_{k}+\eta(x_k-x_{k-1})$, $x_{k}=proj(\hat{x}_{k})$;\\
7: If $\|x_{k}-x_{k-1}\|_1<\epsilon$, break and output $x_k$;\\
8: $k=k+1$, back to step 4.\\
\noindent\rule{\textwidth}{1pt}
\textbf{Remark} 3. HOPMM-II is obtained by adding the ``heavy ball" momentum term $\eta(x_k-x_{k-1})$ to higher-order power method. In the HOPMM-II, we will execute the momentum extrapolation at every 2 steps.

How to choose the parameters $\beta,\eta$ is crucial for the performance of HOPMM-I and HOPMM-II. However, it is difficult to select the parameter $\beta,\eta$ so far. Thus we further propose a free-parameter quadratic extrapolation method for solving the tensor equation (\ref{approximated model}) in the following subsection.
\subsection{Higher-order quadratic extrapolation method}
In this subsection, we extend the quadratic extrapolation method in \cite{Kamvar03} for solving the tensor equation (\ref{approximated model}), referred to as the QEHOPM.
For the classic quadratic extrapolation method, Sidia in \cite{Sidia08} has proved that this method is faster than power method.
We develop the Quadratic Extrapolation Higher Order Power Method (QEHOPM) as follows.

Letting $B=\mathcal{P}\bar{x}^{m-2}$, where $\bar{x}$ is a solution of equation (\ref{approximated model}),
 it is clear that $B$ is also a Markov Matrix, and $\bar{x}$ is the principal eigenvector of $B$.
Assume that the matrix $B=\mathcal{P}\bar{x}^{m-2}$ has only 3 eigenvectors.
Then, the iterate $x_{k-3}$
can be expressed as a linear combination of
these 3 eigenvectors. Of course, $B$
has more than 3 eigenvectors, and $x_{k-3}$
can only be approximated as a linear combination of the first three eigenvectors.
%
 Therefore, the
$\hat{x}$ that we compute in QEHOPM algorithm is only an
estimate for the true
$\bar{x}$.

Similar to the quadratic extrapolation method in \cite{Kamvar03}, we assume that $B$
has only three eigenvectors,
and approximating $x_{k-3}$
as a linear combination of these three eigenvector.
We then define the successive iterates
$x_{k-2}=Bx_{k-3}$, $x_{k-1}=Bx_{k-2}$, $x_{k}=Bx_{k-1}$. From the analysis in Theorem \ref{mainTheorem2}, the QEHOPM will converge to the fixed-point $\bar{x}$. So, in practice, we can use $\mathcal{P}x_{k}^{m-2}$ to approximate $\mathcal{P}\bar{x}^{m-2}$ in $k_{th}$ iteration.

Since we assume $B$
has 3 eigenvectors, the characteristic polynomial $p_B(\lambda)$ is given by
$$p_B(\lambda)=\gamma_0+\gamma_1\lambda+\gamma_2\lambda^2+\gamma_3\lambda^3$$

Moreover, since $\lambda=1$ is the leading eigenvalue,
\begin{equation}\label{eq:QcharcteristicPolynomial}
p_B(\lambda)=(\lambda-1)(\beta_0+\beta_1\lambda+\beta_2 \lambda^2)\triangleq (\lambda-1)q_B(\lambda),
\end{equation}
where $\beta_0=\gamma_1+\gamma_2+\gamma_3$, \ $\beta_1=\gamma_2+\gamma_3$, \ $\beta_2=\gamma_3$.

By the Cayley-Hamilton Theorem, for any vector $z\in \Re^n$, we have
$$p_B(B)z=(\gamma_0I+\gamma_1B+\gamma_2B^2+\gamma_3B^3)z=0.$$

Letting $z=x_{k-3}$, we obtain that

\begin{equation}\label{eq:linearcombination}
\gamma_0x_{k-3}+\gamma_1x_{k-2}+\gamma_2x_{k-1}+\gamma_3x_{k}=0.
\end{equation}

Since $p_B(1)=\gamma_0+\gamma_1+\gamma_2+\gamma_3=0$, we have $\gamma_0=-(\gamma_1+\gamma_2+\gamma_3)$.
Letting $y_k=x_k-x_{k-3},\ y_{k-1}=x_{k-1}-x_{k-3},\ y_{k-2}=x_{k-2}-x_{k-3}$, combing with (\ref{eq:linearcombination}), we have
$$\gamma_1y_{k-2}+\gamma_2y_{k-1}+\gamma_3y_{k}=0.$$
Fixing $\gamma_3=1$, and then
\begin{equation}\label{eq:overestmationequation}
\gamma_1y_{k-2}+\gamma_2y_{k-1}=-y_{k}.
\end{equation}
By Least-Square method and QR factorization, we can compute approximatively the above overdetermined system.\\

Again, by the Cayley-Hamilton Theorem, for any vector $z\in \Re^n$, follows from (\ref{eq:QcharcteristicPolynomial}), $q_B(B)z$ is
the eigenvector of $B$ corresponding to eigenvalue 1 (the
principal eigenvector). Letting $z=x_{k-2}$, we have
\begin{equation}
q_B(B)x_{k-2}=(\beta_0I+\beta_1B+\beta_2 B^2)x_{k-2}=\beta_0x_{k-2}+\beta_1x_{k-1}+\beta_2 x_{k}.
\end{equation}
Scaling the above equation by $1/(\beta_0+\beta_1+\beta_2)$, we have
\begin{equation}\label{eq:2monomtum}
\begin{aligned}
\bar{x}&\approx \frac{\beta_0}{\beta_0+\beta_1+\beta_2}x_{k-2}+\frac{\beta_1}{\beta_0+\beta_1+\beta_2}x_{k-1}+\frac{\beta_2}{\beta_0+\beta_1+\beta_2}x_{k}\\
&=x_{k}+\frac{\beta_0} {\beta_0+\beta_1+\beta_2}(x_{k-2}-x_{k-1})+\frac{\beta_0+\beta_1}{\beta_0+\beta_1+\beta_2}(x_{k-1}-x_k).
\end{aligned}
\end{equation}
So, according to (\ref{eq:2monomtum}), we know that QEHOPM algorithm is the power method with 2 momentum terms. In particular, all of the parameters could be calculated in closed form.

Now, the QEHOPM algorithm is shown as follows.

\noindent\rule{\textwidth}{1pt}
\underline{\textbf{QEHOPM}} \\
1. Given a transition probability tensor $\mathcal{P}$, maximum $k_{max}$, $\beta>0$, termination tolerance $\epsilon$ and an initial point $x_0$;\\
2. Initialize $k=1$.\\
3. \textbf{repeat}\\
4. \ \ \ $x_{k}=\mathcal{P}x_{k-1}^{m-1}$;\\
5. \ \ \ $\delta=\|x_{k}-x_{k-1}\|$;\\
6. \ \ \ periodically,\\
7. \ \ \ \ \ \  $\hat{x}_{k}=Quadratic\ Extrapolation(x_{k-3},\ldots,x_{k})$,\ $x_k=proj(\hat{x}_{k})$;\\
8.\ \ \ \ $k=k+1$;\\
9. \textbf{until} $\delta<\epsilon$\\
\noindent\rule{\textwidth}{1pt}
The quadratic extrapolation algorithm is defined as follows.

\noindent\rule{\textwidth}{1pt}\\
\underline{\textbf{Quadratic Extrapolation}} \\
function $\hat{x}=Quadratic\ Extrapolation(x_{k-3},\ldots,x_{k})\{$\\
for $j=k-2:k \ do$\\
\ \ $y_j=x_j-x_{k-3}$;\\
end\\
$Y=(y_{k-2}\ \ y_{k-1})$;
$\gamma_3=1$;\\
 $\left(\begin{array}{ccc}
\gamma_1\\\gamma_2
\end{array}\right)=-Y^+y_{k}$;\\
 $\beta_0=\gamma_1+\gamma_2+\gamma_3$;\\
$\beta_1=\gamma_2+\gamma_3$;\\ $\beta_2=\gamma_3$;\\
$\hat{x}=\frac{\beta_0}{\beta_0+\beta_1+\beta_2}x_{k-2}+\frac{\beta_1}{\beta_0+\beta_1+\beta_2}x_{k-1}+\frac{\beta_2}{\beta_0+\beta_1+\beta_2}x_{k};$\\
\}\\
\noindent\rule{\textwidth}{1pt}
Using the following Gram-Schmidt to solve $\gamma_1$ and $\gamma_2$.\\
\noindent\rule{\textwidth}{1pt}
\underline{\textbf{Gram-Schmidt}} \\
1. Compute the reduced $QR$ factorization $Y=QR$ using 2 steps of Gram-Schmidt.\\
2. Compute the vector $-Q^Ty_{k}$.\\
3. Solve the upper triangular system:\\
$
R \left(\begin{array}{ccc}
\gamma_1\\\gamma_2
\end{array}\right)=-Q^Ty_{k}
$;\\
for
$R \left(\begin{array}{ccc}
\gamma_1\\\gamma_2
\end{array}\right)$ using back substitution.\\
\noindent\rule{\textwidth}{1pt}\\

In this paper, we will apply quadratic extrapolation at every 4 steps.
\section{Convergence analysis for the proposed methods}
In this section, we present the convergence analysis of the proposed algorithms. Before giving these Theorems, some lemmas that established by Li et al. in \cite{WMNG14,Liu19} are shown as follows.
\begin{lemma}\label{dl2}
If $\mathcal{P}$ is a non-negative transition probability tensor of order $m$ and dimension $n$, then there exists a nonzero non-negative vector
$\bar{x}$ satisfies (\ref{approximated model}). In particular, if $\mathcal{P}$ is irreducible, then $\bar{x}$ must be positive.
\end{lemma}
\begin{lemma}\label{YL2}
Suppose $\mathcal{P}$ is a non-negative transition probability tensor of order $m$ and dimension $n$. If $\delta_m>\frac{m-2}{m-1}$, the $\delta_m$ is given as follows
\begin{equation}\label{gamma3}
\delta_m :=\min_{S\subset\langle n\rangle}\left\{\min_{i_2,\ldots,i_m\in\langle n\rangle}\sum_{i\in S'}p_{i,i_2,\ldots,i_m}+\min_{i_2,\ldots,i_m\in\langle n\rangle}\sum_{i\in S}p_{i,i_2,\ldots,i_m}\right\}.
\end{equation}
where $\langle n\rangle=\{1,2,\ldots,n\}$, $S$ is a subset of $\langle n\rangle$ and $S'$ be its complementary set in $\{1,2,\ldots,n\}$, i,e., $S'=\{1,2,\ldots,n\}\backslash S$.
then the nonzero non-negative vector $\bar{x}$ in Lemma \ref{dl2} is unique.
\end{lemma}

Lemma \ref{dl2} and Lemma \ref{YL2} give the existence and uniqueness conditions of the solution for equation (\ref{approximated model}), respectively.
\begin{lemma}\label{dl4}
Suppose $\mathcal{P}$ is a non-negative transition probability tensor of order $m$ and dimension $n$ and $x,y\in R^n$ are transition probability vectors. Then we have
\begin{equation}\label{inequatiltiy}
\|\mathcal{P}(x^{m-1}-y^{m-1})\|_1\leq \eta_m\|x-y\|_1,
\end{equation}
where $\eta_m=(1-\delta_m)(m-1)$.
\end{lemma}
The proof of Lemma \ref{dl4} can be found in the Lemma 2 of \cite{Liu19}.
\begin{lemma}\label{dl5}
Let $\hat{x},y\in R^n$ and $\|\hat{x}\|_1=1,\|y\|_1=1$. If $x=proj(\hat{x})$, then $\|\hat{x}-y\|_1\geq\|x-y\|_1.$
\end{lemma}
The proof can be obtained by Lemma 3 of \cite{Liu19}.

Based on these above Lemmas, we establish the following convergence Theorems for HOPMM-I, HOPMM-II and QEHOPM, respectively.
\subsection{Convergence analysis for HOPMM-I and HOPMM-II}

\begin{theorem}\label{mainTheorem 1}
Let $\mathcal{P}$ be a non-negative transition probability tensor of order $m$ and dimension $n$ with $\delta_m>\frac{m-2}{m-1}$ and $\bar{x}$ is a solution of equation (\ref{approximated model}). Then, if $\beta<1-\eta_m$ the iterative sequence $\{x_k\}$ generated by HOPMM-I exists a convergent subsequence $\{x_{k_n}\}$ that converges to the solution $\bar{x}$ for any initial transition probability vector $x_0$, i.e.,
\begin{equation}\label{convergence 1}
\lim_{n\rightarrow\infty}x_{k_n}=\bar{x}.
\end{equation}
\end{theorem}
\textit{Proof}. According to the Lemma \ref{YL2} and condition $\delta_m>\frac{m-2}{m-1}$, we get that equation (\ref{approximated model}) has a unique solution $\bar{x}$. From the HOPMM-I, it is easy to get that $x_k\geq0$ for all $k$.

Let $\hat{e}_{k}=\hat{x}_{k}-\bar{x}$ and $e_{k}=x_{k}-\bar{x}$. By Algorithm 1, we can obtain
\begin{equation}\label{Alagequation}
\begin{aligned}
        \hat{e}_{k} &= x_{k}+\beta (x_{k-1}-x_{x-2})-\bar{x},\\
    \end{aligned}
\end{equation}
where $\beta>0$.

By substituting the $\bar{x}=P\bar{x}^{m-1}$ and $x_k=\mathcal{P}x_{k-1}^{m-1}$ into (\ref{Alagequation}), we have
\begin{equation}\label{equation44}
\begin{aligned}
        \hat{e}_{k} &=\mathcal{P}(x_{k-1}^{m-1}-\bar{x}^{m-1})+\beta(x_{k-1}-\bar{x})+\beta(\bar{x}-x_{k-2})\\
        &=\mathcal{P}(x_{k-1}^{m-1}-\bar{x}^{m-1})+\beta e_{k-1}-\beta e_{k-2}.\\
    \end{aligned}
\end{equation}

By Lemma \ref{dl4}, we have
\begin{equation}\label{equation33}
 \|\mathcal{P}(x_{k-1}^{m-1}-\bar{x}^{m-1})\|_1\leq \eta_m\|e_{k-1}\|_1.
\end{equation}
and
\begin{equation}\label{equation444}
  \|e_{k-1}\|_1=\|\mathcal{P}x_{k-2}^{m-1}-\mathcal{P}\bar{x}^{m-1}\|_1\leq\eta_m\|e_{k-2}\|_1
\end{equation}

Then, by (\ref{equation33}) and (\ref{equation444}), we have
\begin{equation}\label{equation55}
\begin{aligned}
\|\hat{e}_{k}\|_1&\leq(\eta_m+\beta)\|e_{k-1}\|_1+\beta\|e_{k-2}\|_1\\
        &\leq(\eta_m+\beta)\eta_m\|e_{k-2}\|_1+\beta\|e_{k-2}\|_1\\
\end{aligned}
\end{equation}

By Lemma (\ref{dl5}), we can obtain
\begin{equation}\label{equation66}
  \|e_{k}\|_1\leq\|\hat{e}_{k}\|_1\leq[(\eta_m+\beta)\eta_m+\beta]\|e_{k-2}\|_1.
\end{equation}

If $\beta<1-\eta_m$, it is easy to get $(\eta_m+\beta)\eta_m+\beta<1$, which proves that the the iterative sequence $\{x_k\}$ generated by HOPMM exists a convergent subsequence $\{x_{k_n}\}$.

\begin{theorem}\label{mainTheorem 22}
Let $\mathcal{P}$ be a non-negative transition probability tensor of order $m$ and dimension $n$ with $\delta_m>\frac{m-2}{m-1}$ and $\bar{x}$ is a solution of equation (\ref{approximated model}). Then, if $\eta<\frac{1-\delta_m}{1+\delta_m}$, the iterative sequence $\{x_k\}$ generated by HOPMM-II converges to the solution $\bar{x}$ for any initial transition probability vector $x_0$, Furthermore, we have the following error bound
\begin{equation}\label{convergence 211}
\|x_{k}-\bar{x}\|_1\leq\epsilon_\eta^k\|x_0-\bar{x}\|_1,
\end{equation}
where $\epsilon_\eta=(\eta+1)\eta_m+\eta$.
\end{theorem}
\textit{Proof}. According to the Theorem \ref{mainTheorem 1}, we get that equation (\ref{approximated model}) has a unique solution $\bar{x}$. From the HOPMM-II, it is obvious that $x_k\geq0$ for all $k$.

Let $\hat{e}_{k}=\hat{x}_{k}-\bar{x}$ and $e_{k}=x_{k}-\bar{x}$. By Algorithm 2, we obtain
\begin{equation}\label{Alagequation 2}
\begin{aligned}
        \hat{e}_{k} &= x_{k}+\eta(x_{k}-x_{x-1})-\bar{x},\\
    \end{aligned}
\end{equation}
where $\eta>0$.

By substituting the $\bar{x}=P\bar{x}^{m-1}$ and $x_k=\mathcal{P}x_{k-1}^{m-1}$ into (\ref{Alagequation 2}), we have
\begin{equation}\label{equation244}
        \hat{e}_{k} =\mathcal{P}(x_{k-1}^{m-1}-\bar{x}^{m-1})+\eta\mathcal{P}(x_{k-1}^{m-1}-\bar{x}^{m-1})-\eta e_{k-1}
\end{equation}

By Lemma \ref{dl4}, we have
\begin{equation}\label{equation233}
 \|\mathcal{P}(x_{k-1}^{m-1}-\bar{x}^{m-1})\|_1\leq \eta_m\|e_{k-1}\|_1.
\end{equation}

Then, by (\ref{equation244}) and (\ref{equation233}), we have
\begin{equation}\label{equation255}
\|\hat{e}_{k}\|_1\leq[(1+\eta)\eta_m+\eta]\|e_{k-1}\|_1£¬
\end{equation}

By Lemma (\ref{dl5}), we get
\begin{equation}\label{equation266}
  \|e_{k}\|_1\leq\|\hat{e}_{k}\|_1\leq[(1+\eta)\eta_m+\eta]\|e_{k-1}\|_1.
\end{equation}

It follows from (\ref{equation266}) that the error bound (\ref{convergence 211}) holds. It is obvious to get $0<\epsilon_\eta<1$ if $0<\eta<\frac{1-\eta_m}{1+\eta_m}$, which proves the convergence theorem. This completes the proof of the theorem.\ \ \ \ \ \ \ \ \ \ \ \ \ \ \ \ \ \ \ \ \ \ \ \ \ \ \ \ \ \ \ \ \ \ \ \ \ \ \ \ \ \ \ \ \ \ \ \ \ \ \ \ \ \ \ \ \ \ \ \ \ \ \ \ \ \ \ \ \ \ \ \ \ \ \ \ \ \ \ \ \ \ \ \ \ \ \ \ \ \ \ \ \ \ \ \ \ \ \ \ \  $\Box$

\subsection{Convergence analysis for QEHOPM}

Now, we establish the convergence Theorem for QEHOPM.
\begin{theorem}\label{mainTheorem2}
Assume $\mathcal{P}$ is a non-negative transition probability tensor of order $m$ and dimension $n$ with $\delta_m>\frac{m-2}{m-1}$ and $\bar{x}$ is a solution of equation (\ref{approximated model}). Then, the iterative sequence $\{x_k\}$ generated by QEHOPM has a convergent subsequence $\{x_{k_n}\}$ that converges to the solution $\bar{x}$ for any initial transition probability vector $x_0$, i.e.,
\begin{equation}\label{convergence}
\lim_{n\rightarrow\infty}x_{k_n}=\bar{x}.
\end{equation}
\end{theorem}
\textit{Proof}. According to the Lemma \ref{YL2} and condition $\delta_m>\frac{m-2}{m-1}$, we get that equation (\ref{approximated model}) has a unique solution $\bar{x}$. From the QEHOPM, it is easy to get that $x_k\geq0$ for all $k$.

Let $\hat{e}_k=\hat{x}_k-x$ and $e_k=x_k-x$. By Algorithm 2, we can obtain
\begin{equation}\label{equation1}
\begin{aligned}
        \hat{e}_k &= \alpha_1x_k+\alpha_2x_{k-1}+\alpha_3x_{k-2}-\bar{x}\\
    \end{aligned}
\end{equation}
where $\alpha_1=\frac{\beta_2}{\beta_0+\beta_1+\beta_2},\alpha_2=\frac{\beta_1}{\beta_0+\beta_1+\beta_2},\alpha_3=\frac{\beta_0}{\beta_0+\beta_1+\beta_2}$.

By substituting the $x_k=\mathcal{P}x_{k-1}^{m-1},\bar{x}=P\bar{x}^{m-1}$ into (7), we have
\begin{equation}\label{equation2}
\begin{aligned}
        \hat{e}_k &= \alpha_1\mathcal{P}x_{k-1}^{m-1}+\alpha_1P\bar{x}^{m-1}-\alpha_1P\bar{x}^{m-1}+\alpha_2x_{k-1}+\alpha_3x_{k-2}-\bar{x}\\
        &=\alpha_1\mathcal{P}(x_{k-1}^{m-1}-\bar{x}^{m-1})+(\alpha_1-1)\bar{x}+\alpha_2x_{k-1}+(1-\alpha_1-\alpha_2)x_{k-2}\\
        &=\alpha_1\mathcal{P}(x_{k-1}^{m-1}-\bar{x}^{m-1})+(1-\alpha_1)(x_{k-2}-\bar{x})+\alpha_2(x_{k-1}-\bar{x}+\bar{x}-x_{k-2})\\
        &=\alpha_1\mathcal{P}(x_{k-1}^{m-1}-\bar{x}^{m-1})+(1-\alpha_1-\alpha_2)e_{k-2}+\alpha_2e_{k-1}\\
    \end{aligned}
\end{equation}

By Lemma \ref{dl4}, we have
\begin{equation}\label{equation3}
 \|\mathcal{P}(x_{k-1}^{m-1}-\bar{x}^{m-1})\|_1\leq \eta_m\|e_{k-1}\|_1.
\end{equation}
and
\begin{equation}\label{equation4}
  \|e_{k-1}\|_1=\|\mathcal{P}x_{k-2}^{m-1}-\mathcal{P}\bar{x}^{m-1}\|_1\leq\eta_m\|e_{k-2}\|_1
\end{equation}

Then, by (\ref{equation2}),(\ref{equation3}) and (\ref{equation4}), we have
\begin{equation}\label{equation5}
\begin{aligned}
\|\hat{e}_k\|_1&\leq(\alpha_1\eta_m+\alpha_2)\|e_{k-1}\|_1+\alpha_3\|e_{k-2}\|_1\\
&=[(\alpha_1\eta_m+\alpha_2)\eta_m+\alpha_3]\|e_{k-2}\|_1.\\
\end{aligned}
\end{equation}

By Lemma (\ref{dl5}), we can obtain
\begin{equation}\label{equation6}
  \|e_{k}\|_1\leq\|\hat{e}_k\|_1\leq[(\alpha_1\eta_m+\alpha_2)\eta_m+\alpha_3]\|e_{k-2}\|_1.
\end{equation}

It is easy to get that $0\leq\eta_m<1$ when $\delta_m>\frac{m-2}{m-1}$. Then, we have
\begin{equation}\label{equation7}
  (\alpha_1\eta_m+\alpha_2)\eta_m+\alpha_3<1.
\end{equation}

Now, by (\ref{equation6}) and (\ref{equation7}), we can get that sequence $\{e_k\}$ has a convergent subsequence $\{e_{k_n}\}$ that will converges to zero vector, which proves that
the iterative sequence $\{x_k\}$ has a convergent subsequence $\{x_{k_n}\}$ that converges to the solution $\bar{x}$. This completes the proof of the theorem.\ \ \ \ \ \ \ \ \ \ \ \ \ \ \ \ \ \ \ \ \ \ \ \ \ \ \ \ \ \ \ \ \ \ \ \ \ \ \ \ \ \ \ \ \ \ \ \ \ \ \ \ \ \ \ \ \ \ \ \ \ \ \ \ \ \ \ \ \ \ \ \ \ \ \ \ \ \ \ \ \ \ \ \ \ \ \ \ \ \ \ \ $\Box$
\section{Numerical experiments}
In this section, a number of numerical experiments are presented to verify the efficiency and superiority of our methods, compared with
the original higher-order power method(HOPM)\cite{Liu19}, the relaxation higher-order power method((RHOPM))\cite{Liu19}, the shifted power method(S)\cite{Gleich14} and the inner-outer iteration method(IO)\cite{Gleich14}. Three measure indexes are reported, including the number of iterations(denoted IT), the CPU time in seconds(denoted by CPU) and the relative residual(denoted by RR) defined by
$\|\mathcal{P}x_{k}^{m-1}-x_{k}\|_1.$
\par In the numerical experiments, all initial points are chosen to be $x_0=ones(n,1)/n$, all algorithms are performed with Tensor Toolbox 2.6 \cite{BaderKolda2015} in MATLAB R2010a and are terminated when the condition
$\|x_{k+1}-x_k\|<10^{-10}$
is satisfied. The maximum iterative number is set to 1000. The curve of the norm of relative residual vector versus the number of iteration step is plotted. The selection of parameter in RHOPM are the same to that of in \cite{Liu19}.

\subsection{Numerical results for higher-order Markov chains}
In this subsection, we use the proposed methods(i.e., QEHOPM, HOPMM-I and HOPMM-II), HOPM, RHOPM and GEAP for solving the limiting probability distribution vector of four transition probability tensors (which were contained in the appendix).

\begin{figure}[H]
$$
\begin{array}{cccc}
\includegraphics[width=0.45\textwidth]{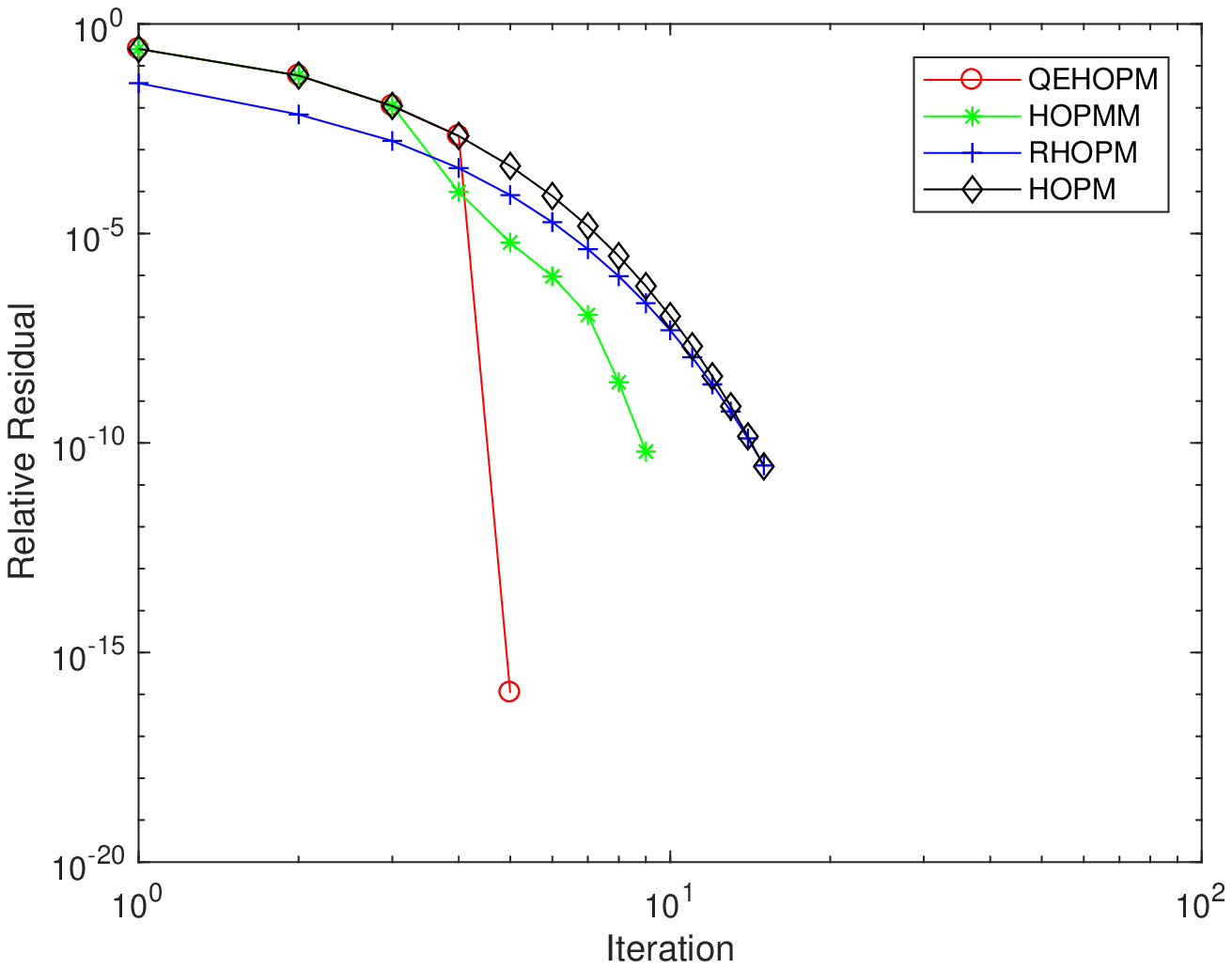}&
\includegraphics[width=0.45\textwidth]{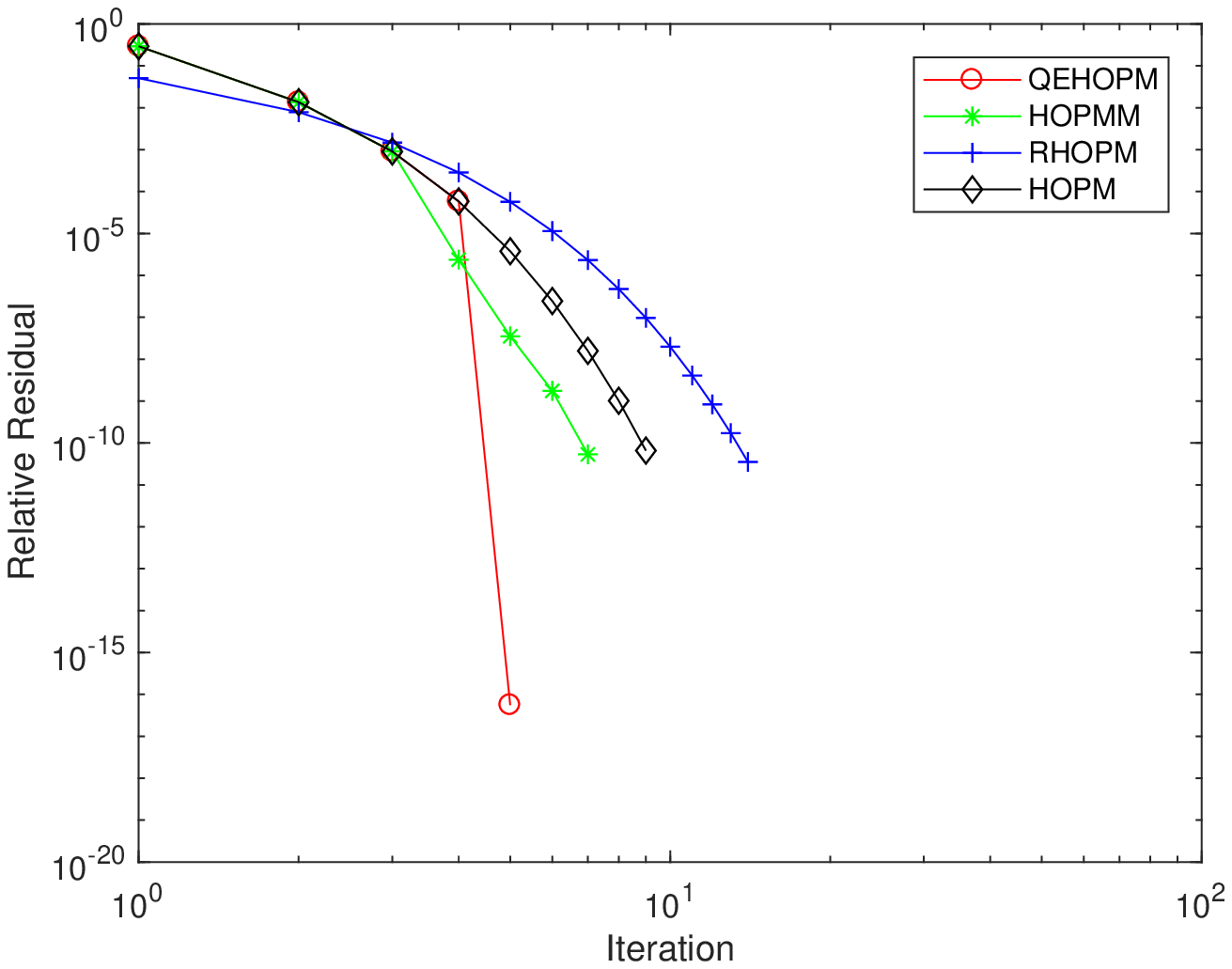}\\
(i) &  (ii) \\
\includegraphics[width=.45\textwidth]{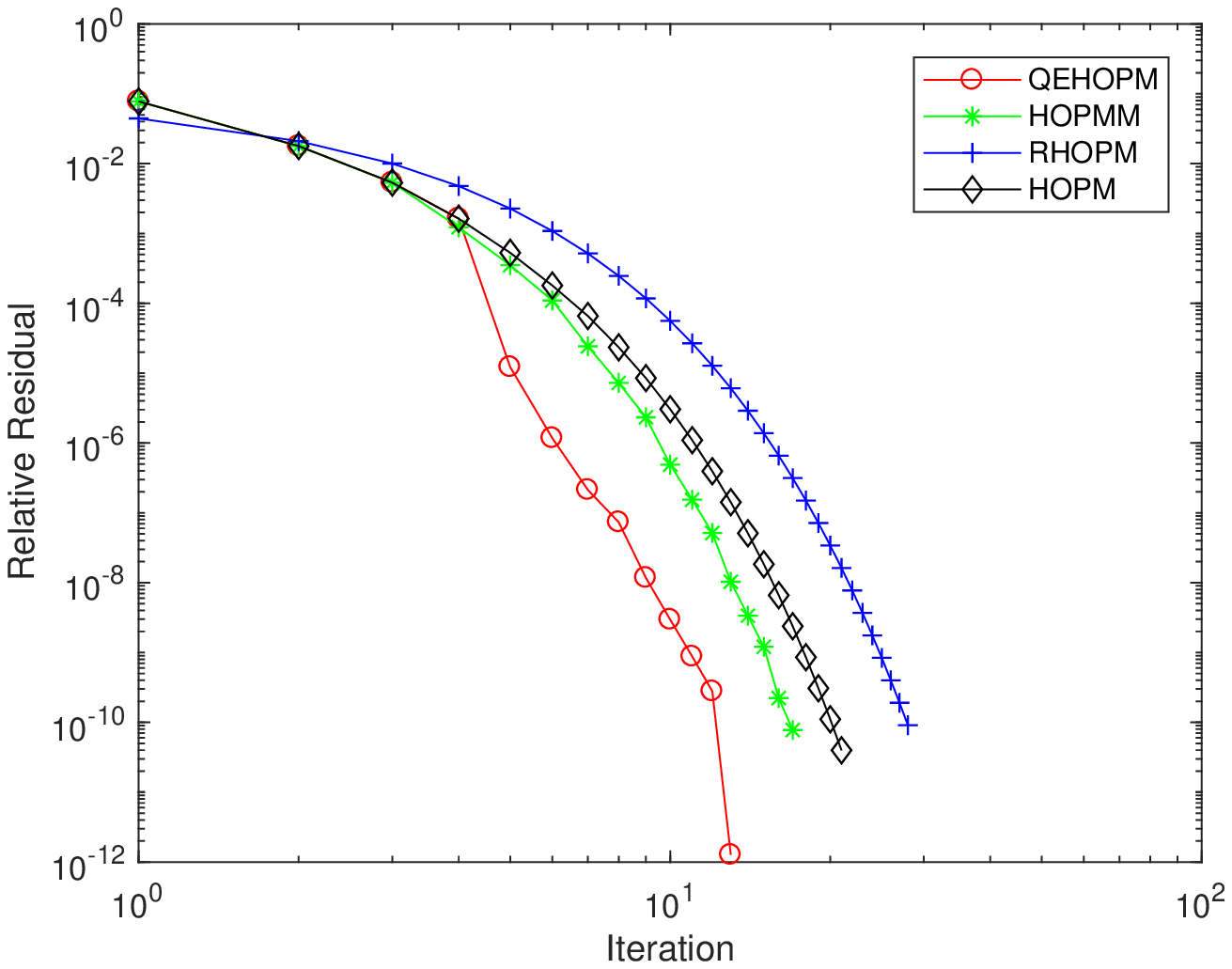}&
\includegraphics[width=.45\textwidth]{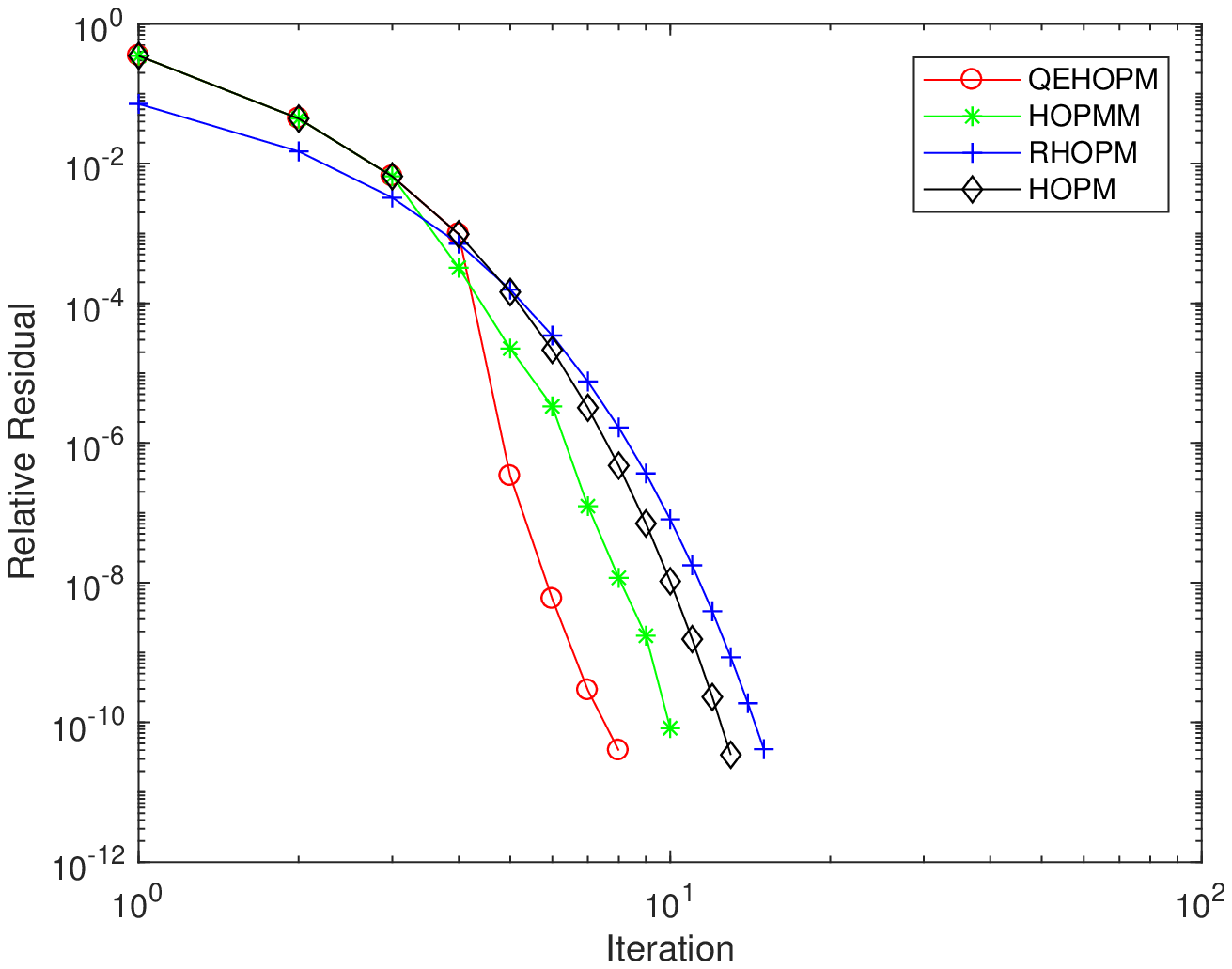}\\
(iii)&(iv)
\end{array}
$$
\caption{ The curve of the norm of relative residual vector versus the number of iteration step for (i)-(iv)}.
\end{figure}

\begin{table}[H]
\centering
\caption{The computed quantities of the examples of tensors (i)-(iv).}
\begin{tabular}{c|cccccc}
\hline
Examples                  & Algorithm           & CPU          & IT     & RR\\
\hline
(i)                       &HOPM         & 0.0148        & 15   &  2.76e-11  \\
                          &GEAP         & 0.0579        & 14   &  9.06e-11  \\
                          &RHOPM($\gamma=1.2$)         & 0.0170        & 16   &  2.90e-11  \\
                          &HOPMM-I   ($\beta=0.045$)              & 0.0105        & 9    & 6.21e-11 \\
                          &HOPMM-II   ($\eta=0.2$)              & 0.0095        & 10    & 3.82e-11 \\
                          &QEHOPM               & 0.0039        &5    & 1.11e-16 \\
\hline
(ii)                      &HOPM         & 0.0079        & 9   &  6.58e-11  \\
                          &GEAP         & 0.0449        & 9   &  7.58e-12  \\
                          &RHOPM($\gamma=1.2$)         & 0.0154        & 15   &  3.52e-11  \\
                          &HOPMM-I ($\beta=0.0045$)               & 0.0066        & 7   & 9.24e-11 \\
                          &HOPMM-II   ($\eta=0.07$)              & 0.0049        & 6    & 3.82e-11 \\
                          &QEHOPM        & 0.0033        &5    & 5.55e-17 \\

\hline
(iii)                     &HOPM         & 0.0168        & 21   &  3.97e-11  \\
                          &GEAP         & 0.1540        & 27   &  5.98e-11  \\
                          &RHOPM($\gamma=1.2$)          & 0.0308        & 29   &  9.08e-11  \\
                          &HOPMM-I ($\beta=0.1$)         & 0.0183        & 17    & 7.76e-11 \\
                          &HOPMM-II   ($\eta=0.1$)       & 0.0189        & 19    & 3.82e-11 \\
                          &QEHOPM              & 0.0071        &13    & 1.26e-12 \\

\hline
(iv)                      &HOPM         & 0.0147       & 13   &  3.42e-11  \\
                          &GEAP         & 0.0949       &12   &  9.52e-11  \\
                          &RHOPM ($\beta=1.2$)        & 0.0189& 16   &  4.12e-11  \\
                          &HOPMM-I ($\beta=0.03$)               & 0.0087        & 10    & 8.25e-11 \\
                          &HOPMM-II   ($\eta=0.2$)              & 0.0098        & 9    & 3.82e-11 \\
                          &QEHOPM                  & 0.0070        &8    & 3.95e-11 \\
\hline
\end{tabular}
\end{table}
\par The numerical results are reported in Table 1. As we can see, from the Table 1, the number of the iteration steps in QEHOPM and HOPMM-I/HOPMM-II are less than that of HOPM, RHOPM and GEAP. Furthermore, our methods(i.e., QEHOPM and HOPMM) spend less time than HOPM, RHOPM and GEAP. Specially, QEHOPM performs the best among all methods.

\par Figure 1 ploted the norm of relative residual vector versus the number of iteration step for the above examples. Compared with HOPM, RHOPM need more iterations while HOPMM-I/HOPMM-II is much faster. Noticed that QEHOPM is the best one. Especially for (i) and (ii), only one quadratic extrapolation can reach to the solution.


\subsection{Numerical results for multilinear PageRank}

In this subsection, we display the numerical results when the tested algorithms are applied for solving the multilinear PageRank. We use the benchmark set of 29 stochastic tensors constructed by Gleich et al\cite{Gleich14}. For the sake of fairness, 
we rewrite the codes of {\bf IO} method (the inner-outer interation in \cite{Gleich14}) and {\bf S} method (the shifted fixed-point iteration in \cite{Gleich14}) by using the function \textit{ttv} of the package Tensor Toolbox 2.6.
The vector $v$ is set to $\frac{1}{n}e$, where $e=ones(n,1)$, and the damping parameter $\theta$ is set to 0.70, 0.85, 0.90, 0.95, and 0.99, respectively. we list the numerical results in Tables 2-6 for different $\theta$, where $'-'$ means that the corresponding algorithm is not available.

\begin{table}[H]
\centering
\caption{The numerical results for multilinear PageRank with $\theta=0.7$.}
\begin{tabular}{cccccc}
\hline
                   &HOPM     &QEHOPM        &RHOPM    &IO      &S \\
                   & IT\ \ CPU  & IT\ \ CPU & IT\ \ CPU  & IT\ \ CPU   & IT\ \ CPU\\
\hline
$R_{3,1}$          & 52\ \ \ 0.1271   &{\bf 9\ \ 0.0183}         &28\ \ \ 0.0468    & 69\ \ \  0.4481    &37\ \ \      0.0645\\
$R_{3,2}$          & 15\ \ \ 0.0949  &{\bf 9\ \  0.0167}         &14\ \ \ 0.0183    & 60\ \ \  0.3419    &30\ \ \      0.0575\\
$R_{3,3}$          & 15\ \ \ 0.1078   &{\bf 9\ \  0.0179}        &14\ \ \ 0.0202    & 60\ \ \  0.3822    &30\ \ \      0.0877\\
$R_{3,4}$          & 43\ \ \ 0.1189   &{\bf 9\ \  0.0101}        &18\ \ \ 0.0551   &48\ \ \   0.3955     &20\ \ \     0.0429\\
$R_{3,5}$          & 86\ \ \ 0.1704   &{\bf 17\ \ 0.0242}        &96\ \ \ 0.1022   &182\ \ \  1.6651    &134\ \ \     0.2289\\
$R_{4,1}$          & 55\ \ \ 0.1317   &{\bf 37\ \ 0.0438}       &62\ \ \ 0.0887     &130\ \ \  1.1107    &88\ \ \      0.1655\\
$R_{4,2}$          & 77\ \ \ 0.1780   &{\bf 34\ \ 0.0405}       &64\ \ \ 0.0866    &128\ \ \  1.0479    &109\ \ \     0.2062\\
$R_{4,3}$         & 50\ \ \ 0.1300   &{\bf 37\ \  0.0678}       &51\ \ \ 0.0777   &114\ \ \  0.9529    &79\ \ \      0.1481 \\
$R_{4,4}$         & 77\ \ \ 0.1521    &{\bf59\ \   0.0563}     &65\ \ \ 0.0846     &118\ \    1.0502     &80\ \ \     0.1689\\
$R_{4,5}$          & 53\ \ \ 0.1254   &{\bf 34\ \  0.0386}      &40\ \ \ 0.0648     &115\ \ \  0.9562    &81\ \ \      0.1530\\
$R_{4,6}$         & 64\ \ \ 0.1417   &{\bf 33\ \  0.0352}       &51\ \ \ 0.0776    &117\ \ \  0.8601    &93\ \ \      0.1785\\
$R_{4,7}$         & 51\ \ \ 0.1279   &{\bf 44\ \  0.0450}       &53\ \ \ 0.0799    &109\ \ \  0.9250    &75\ \ \      0.1347\\
$R_{4,8}$         & 68\ \ \ 0.1421   &{\bf 56\ \  0.0570}       &61\ \ \ 0.0890    &112\ \ \  0.8877    &94\ \ \      0.1534\\
$R_{4,9}$         & 52\ \ \ 0.1359   &{\bf 27\ \  0.0406}       &45\ \ \ 0.0628    &110\ \ \  0.9428    &81\ \ \      0.1365\\
$R_{4,10}$         & 90\ \ \ 0.2322   & 67\ \ 0.0904       &{\bf57\ \ \ 0.0742}    &111\ \ \  0.9534    &80\ \ \      0.1535\\
$R_{4,11}$         & 54\ \ \ 0.1327   &{\bf 38\ \  0.0454}      &49\ \ \ 0.0715    & 114\ \ \ 0.8166   &82\ \ \      0.1408\\
$R_{4,12}$         & 38\ \ \ 0.1215   &{\bf 31\ \  0.0400}      &39\ \ \ 0.0628     & 86\ \ \  0.8176    &61\ \ \      0.1156\\
$R_{4,13}$         & 58\ \ \ 0.1333   &{\bf 49\ \ 0.0488}       &64\ \ \ 0.0800    & 125\ \ \  1.0120   &86\ \ \      0.1468 \\
$R_{4,14}$         & 56\ \ \ 0.1175   &{\bf 20\ \ 0.0269}       &41\ \ \ 0.0641    & 111\ \ \  0.8955   &82\ \ \      0.1403\\
$R_{4,15}$         & 70\ \ \ 0.1494   &{\bf 37\ \  0.0422}      &59\ \ \ 0.0788    & 122\ \ \  0.9808   &101\ \ \     0.1638\\
$R_{4,16}$         & 56\ \ \ 0.1342   &{\bf 49\ \  0.0496}      &62\ \ \ 0.0799    & 118\ \ \  0.8728   &85\ \ \      0.1481\\
$R_{4,17}$         & 53\ \ \ 0.1209   &{\bf 23\ \   0.0351}     &56\ \ \ 0.0709    & 114\ \ \  0.8236   &78\ \ \      0.1273\\
$R_{4,18}$          & 52\ \ \ 0.1323   &{\bf 29\ \   0.0316}    &55\ \ \ 0.0719    & 116\ \ \  0.8247   &78\ \ \      0.1386\\
$R_{4,19}$        & 49\ \ \ 0.1277   &{\bf 40\ \   0.0436}      &54\ \ \ 0.0722    & 103\ \ \  0.8040   &75\ \ \      0.1312\\
$R_{6,1}$          & 40\ \ \ 0.1124   & 26\ \ {\bf 0.0266 }     &{\bf24}\ \ \ 0.0442     & 100\ \ \  0.6884   &65\ \ \      0.1045\\
$R_{6,2}$          & 49\ \ \ 0.1234   &{\bf 25\ \  0.0290}      &44\ \ \ 0.0773     & 110\ \ \  0.7558   &64\ \ \      0.1175\\
$R_{6,3}$          & 35\ \ \ 0.1091   &{\bf 22\ \   0.0304}     &40\ \ \ 0.0667     & 89 \ \ \  0.6202   &57\ \ \      0.0928\\
$R_{6,4}$          & 40\ \ \ 0.1180   &{\bf 21\ \  0.0289}      &25\ \ \ 0.0404     & 104\ \ \  0.7710   &66\ \ \      0.1096\\
$R_{6,5}$          & 33\ \ \ 0.1033   &{\bf 19\ \  0.0270}      &26\ \ \ 0.0455     & 87\ \ \  0.4741    &55\ \ \      0.0946\\
\hline
\end{tabular}
\end{table}

\begin{table}[H]
\centering
\caption{The numerical results for multilinear PageRank with $\theta=0.85$.}
\begin{tabular}{cccccc}
\hline
                   &HOPM     &QEHOPM        &RHOPM    &IO      &S \\
                   & IT\ \ CPU  & IT\ \ CPU & IT\ \ CPU  & IT\ \ CPU   & IT\ \ CPU\\
\hline
$R_{3,1}$          & 130\ \ \ 0.2141   &{\bf 13\ \ 0.0193}         &34\ \ \ 0.0368     & 60\ \ \  0.5886      &44\ \ \      0.0783\\
$R_{3,2}$          & 21\ \ \ 0.0843   &{\bf 13\ \  0.0220}         &29\ \ \ 0.0334     & 55\ \ \  0.3816      &37\ \ \      0.0708\\
$R_{3,3}$          & 21\ \ \ 0.0822   &{\bf 13\ \  0.0218}         &17\ \ \ 0.0220     & 55\ \ \  0.3814      &37\ \ \      0.0681\\
$R_{3,4}$          & 70\ \ \ 0.1359   &{\bf 13\ \  0.0203}         &19\ \ \ 0.0205     & 39\ \ \  0.3800      &21\ \ \      0.0414\\
$R_{3,5}$          & 103\ \ \ 0.1617   &{\bf 29\ \ 0.0350}        &115\ \ \ 0.1022     & 160\ \ \  1.9054      &158\ \ \      0.2571\\
$R_{4,1}$          & 134\ \ \ 0.1857   &{\bf 102\ \ 0.0978}         &146\ \ \ 0.1287     & 166\ \ \  1.9235      &191\ \ \      0.3081\\
$R_{4,2}$          & 261\ \ \ 0.2957   &{\bf 33\ \ 0.0356}         &318\ \ \ 0.2866     & 158\ \ \  1.8214      &296\ \ \      0.4762\\
$R_{4,3}$          & 101\ \ \ 0.1550   &{\bf 70\ \  0.0702}         &125\ \ \ 0.1137     & 147\ \ \  1.4540      &155\ \ \      0.2752 \\
$R_{4,4}$          & 262\ \ \ 0.3049   &172\ \  {\bf 0.1639}         &181\ \ \ 0.1646     & {\bf 167}\ \   1.9003      &232\ \ \      0.3767\\
$R_{4,5}$          & 162\ \ \ 0.2413   &{\bf 41\ \  0.0514}         &144\ \ \ 0.1648     & 160\ \ \  2.2219      &220\ \ \      0.4797\\
$R_{4,6}$          & 225\ \ \ 0.2751   &{\bf 29\ \  0.0356}         &225\ \ \ 0.1976     & 155\ \ \  1.6964      &269\ \ \      0.4301\\
$R_{4,7}$          & 120\ \ \ 0.1703   &{\bf 110\ \  0.1044}         &132\ \ \ 0.1159     & 141\ \ \  1.6796      &170\ \ \      0.2723\\
$R_{4,8}$          & 300\ \ \ 0.3332   &{\bf 57\ \  0.0572}         &279\ \ \ 0.2490     & 136\ \ \  1.5656      &296\ \ \      0.4816\\
$R_{4,9}$          & 144\ \ \ 0.1996   &{\bf 49\ \  0.0537}         &132\ \ \ 0.1158     & 144\ \ \  1.6666      &196\ \ \      0.3150\\
$R_{4,10}$          & 364\ \ \ 0.3767   &{\bf 113\ \ 0.1092}         &130\ \ \ 0.1242     & 154\ \ \  1.8034      &172\ \ \      0.2779\\
$R_{4,11}$          & 133\ \ \ 0.2074   &{\bf 64\ \  0.0726}         &115\ \ \ 0.1015     & 174\ \ \  1.6932      &195\ \ \      0.3086 \\
$R_{4,12}$          & 88\ \ \ 0.1495   &{\bf 43\ \   0.0451}         &80\ \ \ 0.0728     & 140\ \ \  1.5340      &137\ \ \      0.2395\\
$R_{4,13}$          & 151\ \ \ 0.2028   &{\bf 122\ \ 0.1161}         &162\ \ \ 0.1428     & 177\ \ \  2.0442      &215\ \ \      0.3434 \\
$R_{4,14}$          & 134\ \ \ 0.1915   &{\bf 33\ \ 0.0363}        &117\ \ \ 0.1041     & 135\ \ \  1.5947      &181\ \ \      0.2915\\
$R_{4,15}$          & 238\ \ \ 0.2808   &{\bf 33\ \  0.0374}         &238\ \ \ 0.2088     & 154\ \ \  1.7142      &276\ \ \      0.4491\\
$R_{4,16}$          & 149\ \ \ 0.1969   &{\bf 101\ \  0.0976}         &138\ \ \ 0.1203     & 167\ \ \  1.9701      &209\ \ \      0.3348\\
$R_{4,17}$          & 208\ \ \ 0.2491   &{\bf 35\ \   0.0377}         &136\ \ \ 0.1209     & 177\ \ \  1.9880      &186\ \ \      0.3062\\
$R_{4,18}$          & 141\ \ \ 0.1897   &{\bf 49\ \   0.0518}         &156\ \ \ 0.1419     & 209\ \ \  2.3378      &215\ \ \      0.3558\\
$R_{4,19}$          & 120\ \ \ 0.1733   &{\bf 92\ \   0.0904}         &113\ \ \ 0.1022     & 125\ \ \  1.4163      &163\ \ \      0.3173\\
$R_{6,1}$          & 106\ \ \ 0.1659   &{\bf 30\ \  0.0347}         &60\ \ \ 0.0542     & 162\ \ \  1.5742      &163\ \ \      0.2793\\
$R_{6,2}$          & 98\ \ \ 0.1679   &{\bf 46\ \  0.0511}         &88\ \ \ 0.0773     & 130\ \ \  1.3974      &123\ \ \      0.2517\\
$R_{6,3}$          & 67\ \ \ 0.1325   &{\bf 28\ \   0.0289}         &74\ \ \ 0.0667     & 109\ \ \  1.0417      &103\ \ \      0.1839\\
$R_{6,4}$          & 85\ \ \ 0.1482   &{\bf 26\ \  0.0293}         &52\ \ \ 0.0504     & 129\ \ \  1.3442      &131\ \ \      0.2406 \\
$R_{6,5}$          & 66\ \ \ 0.1288   &{\bf 26\ \  0.0299}         &47\ \ \ 0.0455     & 113\ \ \  0.9191      &106\ \ \      0.1970\\
\hline
\end{tabular}
\end{table}

\begin{table}[H]
\centering
\caption{The numerical results for multilinear PageRank with $\theta=0.9$.}
\begin{tabular}{cccccc}
\hline
                   &HOPM     &QEHOPM        &RHOPM    &IO      &S \\
                   & IT\ \ CPU  & IT\ \ CPU & IT\ \ CPU  & IT\ \ CPU   & IT\ \ CPU\\
\hline
$R_{3,1}$          & 245\ \ \ 0.2159   &{\bf 13\ \ 0.0105}       &36\ \ \ 0.0268     & 58\ \ \  0.5878    &46\ \ \      0.0849\\
$R_{3,2}$          & 23 \ \ \ 0.0207   &{\bf 13\ \  0.0117}      &17\ \ \ 0.0170     & 54\ \ \  0.4342    &41\ \ \      0.0656\\
$R_{3,3}$          & 23\ \ \ 0.0483    &{\bf 13\ \  0.0127}      &17\ \ \ 0.0191     & 54\ \ \  0.3888    &41\ \ \      0.0754\\
$R_{3,4}$          & 88\ \ \ 0.0752    &{\bf 13\ \  0.0114}      &19\ \ \ 0.0195     &36\ \ \   0.3525    &21\ \ \     0.0369\\
$R_{3,5}$          & 68\ \ \ 0.0540    &{\bf 45\ \ 0.0440}       &75\ \ \ 0.0587     &104\ \ \  1.1205    &106\ \ \     0.1666\\
$R_{4,1}$          & 298\ \ \ 0.2555   &{\bf 117\ \ 0.0934}      &283\ \ \ 0.2287    &216\ \ \  2.8204    &384\ \ \      0.7224\\
$R_{4,2}$          & 906\ \ \ 0.7776   &{\bf 30\ \ 0.0274}      &667\ \ \ 0.5866     &181\ \ \  1.9300    &643\ \ \     0.9784\\
$R_{4,3}$          & 157\ \ \ 0.1373   &{\bf 64\ \  0.0494}     &174\ \ \ 0.1737    &197\ \ \  1.9505    &239\ \ \      0.3715 \\
$R_{4,4}$          & 660\ \ \ 0.5677   &{\bf 53\ \ 0.0452}      &432\ \ \ 0.2646     &213\ \    2.4346    &481\ \ \     0.7326\\
$R_{4,5}$          & 537\ \ \ 0.4451   &{\bf 53\ \  0.0438}     &512\ \ \ 0.4648     &203\ \ \  2.4789    &535\ \ \      0.8161\\
$R_{4,6}$          & 873\ \ \ 0.6930   &{\bf 34\ \  0.0297}     &621\ \ \ 0.5976     &182\ \ \  2.0236    &629\ \ \      1.0454\\
$R_{4,7}$          & 243\ \ \ 0.2177   &{\bf 74\ \  0.0593}      &226\ \ \ 0.1159    &182\ \ \  1.9439    &315\ \ \      0.5120\\
$R_{4,8}$          & -\ \ \ \ \ -         &{\bf 53\ \  0.0433}      &745\ \ \ 0.6490    &152\ \ \  1.8511    &885\ \ \      1.4249\\
$R_{4,9}$          & 329\ \ \ 0.2808   &{\bf 42\ \  0.0353}      &337\ \ \ 0.2705    &173\ \ \  1.9783    &368\ \ \      0.5618\\
$R_{4,10}$         & -\ \ \ \ \ -          &{\bf 108\ \ 0.0876}       &235\ \ \ 0.1642   &202\ \ \  2.3554    &302\ \ \      0.5203\\
$R_{4,11}$         & 233\ \ \ 0.2030   &{\bf 46\ \  0.0389}      &204\ \ \ 0.1615    &251\ \ \ 2.3860    &338\ \ \      0.5273\\
$R_{4,12}$         & 151\ \ \ 0.1294   &{\bf 92\ \  0.0881}      &126\ \ \ 0.1011    &202\ \ \  2.2082    &230\ \ \      0.3806\\
$R_{4,13}$         & 314\ \ \ 0.2661   &{\bf 69\ \  0.0554}       &328\ \ \ 0.2208    &226\ \ \  2.8959   &406\ \ \      0.6183 \\
$R_{4,14}$         & 272\ \ \ 0.1175   &{\bf 20\ \ 0.0269}       &271\ \ \ 0.2241    & 161\ \ \  0.8955   &82\ \ \      0.1403\\
$R_{4,15}$         & 799\ \ \ 0.1494   &{\bf 37\ \ 0.0422}      &599\ \ \ 0.5088    & 178\ \ \  0.9808   &101\ \ \     0.1638\\
$R_{4,16}$         & 299\ \ \ 0.1342   &{\bf 49\ \ 0.0496}      &285\ \ \ 0.1403    & 216\ \ \  0.8728   &85\ \ \      0.1481\\
$R_{4,17}$         & -\ \ \ \ \ -    &{\bf 23\ \  0.0351}      &258\ \ \ 0.1209    & 287\ \ \  0.8236   &78\ \ \      0.1273\\
$R_{4,18}$         & 154\ \ \ 0.1323   &{\bf 29\ \ 0.0316}      &140\ \ \ 0.1009    & 161\ \ \  0.8247   &78\ \ \      0.1386\\
$R_{4,19}$         & 218\ \ \ 0.1277   &{\bf 40\ \ 0.0436}     &213\ \ \ 0.1022    & 141 \ \  0.8040   &75\ \ \      0.1312\\
$R_{6,1}$          & 232\ \ \ 0.1124   &{\bf 26\ \  0.0266}      &123\ \ \ 0.0942     & 300\ \ \  0.6884   &65\ \ \      0.1045\\
$R_{6,2}$          & 118\ \ \ 0.1234   &{\bf 25\ \  0.0290}      &118\ \ \ 0.0773     & 143\ \ \  0.7558   &64\ \ \      0.1175\\
$R_{6,3}$          & 95\ \ \ 0.1091   &{\bf 22\ \   0.0304}     &106\ \ \ 0.0667     & 133\ \ \  0.6202   &57\ \ \      0.0928\\
$R_{6,4}$          & 152\ \ \ 0.1180   &{\bf 21\ \  0.0289}      &96\ \ \ 0.0504     & 164\ \ \  0.7710   &66\ \ \      0.1096\\
$R_{6,5}$          & 98\ \ \ 0.1033   &{\bf 19\ \  0.0270}      &64\ \ \ 0.0455     & 142\ \ \  0.4741    &55\ \ \      0.0946\\
\hline
\end{tabular}
\end{table}

\begin{table}[H]
\centering
\caption{The numerical results for multilinear PageRank with $\theta=0.95$.}
\begin{tabular}{cccccc}
\hline
                   &HOPM     &QEHOPM        &RHOPM    &IO      &S \\
                   & IT\ \ CPU  & IT\ \ CPU & IT\ \ CPU  & IT\ \ CPU   & IT\ \ CPU\\
\hline
$R_{3,1}$          & -\ \ \ \ \ \ -    &{\bf 21\ \ 0.0254}      &45\ \ \ 0.0454     & 56\ \ \  0.4543    & -\ \ \ \ \ \ -  \\
$R_{3,2}$          & 28 \ \ 0.1009     &{\bf 13\ \ 0.0222}      &17\ \ \ 0.0255     & 55\ \ \  0.4514    &46\ \ \      0.0902\\
$R_{3,3}$          & 28\ \ \ 0.1029    &{\bf 13\ \ 0.0224}      &17\ \ \ 0.0262     & 55\ \ \  0.4522    &46\ \ \      0.0900\\
$R_{3,4}$          &114\ \ \ 0.1775    &{\bf 13\ \ 0.0174}      &19\ \ \ 0.0288     &34\ \ \   0.3778    &21\ \ \     0.0410\\
$R_{3,5}$          &43\ \ \ 0.1168     &{\bf 37\ \ 0.0458}      &48\ \ \ 0.1198     & 69\ \ \  0.7245    &69\ \ \     0.1312\\
$R_{4,1}$          & -\ \ \ \ \ \ -    &{\bf 21\ \ 0.0230}      & -\ \ \ \ \ \ -    &331\ \ \  4.4715    & -\ \ \ \ \ \ - \\
$R_{4,2}$          & -\ \ \ \ \ \ -    &{\bf 22\ \ 0.0241}      & -\ \ \ \ \ \ -    &220\ \ \  2.6415    & -\ \ \ \ \ \ - \\
$R_{4,3}$          &410\ \ \ 0.4276    &{\bf 49\ \ 0.0531}      &454\ \ \ 0.4973    &358\ \ \  4.1896    &586\ \ \      0.9188 \\
$R_{4,4}$          & -\ \ \ \ \ \ -    &{\bf 25\ \ 0.0262}      & -\ \ \ \ \ \ -       &323\ \    4.1730    & -\ \ \ \ \ \ - \\
$R_{4,5}$          & -\ \ \ \ \ \ -    &{\bf105\ \ 0.0891}      & -\ \ \ \ \ \ -       &282\ \ \  3.5991    & -\ \ \ \ \ \ -  \\
$R_{4,6}$          & -\ \ \ \ \ \ -    &{\bf 41\ \ 0.0378}      & -\ \ \ \ \ \ -       &224\ \ \  2.6912    & -\ \ \ \ \ \ -   \\
$R_{4,7}$          & -\ \ \ \ \ \ -    &{\bf 37\ \ 0.0350}      & -\ \ \ \ \ \ -       &263\ \ \  3.2073    & -\ \ \ \ \ \ -  \\
$R_{4,8}$          & -\ \ \ \ \ \ -    &{\bf 42\ \ 0.0411}      & -\ \ \ \ \ \ -       &174\ \ \  2.2292    & -\ \ \ \ \ \ -  \\
$R_{4,9}$          & -\ \ \ \ \ \ -    &{\bf 42\ \ 0.0391}      & -\ \ \ \ \ \ -       &246\ \ \  2.9428    & -\ \ \ \ \ \ -  \\
$R_{4,10}$         & -\ \ \ \ \ \ -    &{\bf 32\ \ 0.0315}      & -\ \ \ \ \ \ -       &326\ \ \  4.1776    & -\ \ \ \ \ \ -  \\
$R_{4,11}$         &945\ \ \ 0.8263    &{\bf137\ \ 0.1178}      &873\ \ \ 0.7070    &568\ \ \ 5.7697   & -\ \ \ \ \ \ -         \\
$R_{4,12}$         &209\ \ \ 0.2597    &{\bf 30\ \ 0.0319}      &174\ \ \ 0.1493    &202\ \ \  2.5981    &295\ \ \      0.4853\\
$R_{4,13}$         & -\ \ \ \ \ \ -   &{\bf 43\ \ 0.0402}       & -\ \ \ \ \ \ -    &340\ \ \  4.1508   & -\ \ \ \ \ \ -     \\
$R_{4,14}$         & -\ \ \ \ \ \ -      &{\bf 38\ \ 0.0356}    & -\ \ \ \ \ \ -     &217\ \ \  2.6353   & -\ \ \ \ \ \ -    \\
$R_{4,15}$         & -\ \ \ \ \ \ -      &{\bf 42\ \ 0.0397}    & -\ \ \ \ \ \ -    & 221\ \ \ 3.9453   & -\ \ \ \ \ \ -  \\
$R_{4,16}$         & -\ \ \ \ \ \ -      &{\bf 41\ \ 0.0396}    & -\ \ \ \ \ \ -      & 326\ \ \  0.8728   & -\ \ \ \ \ \ -  \\
$R_{4,17}$         & -\ \ \ \ \ \ -      &{\bf934\ \ 0.7409}    & -\ \ \ \ \ \ -      & -\ \ \ \ \ \ -          & -\ \ \ \ \ \ -  \\
$R_{4,18}$         & 867\ \ \ 0.7706   &{\bf 41\ \ 0.0443}     &812\ \ \ 0.7020     &280\ \ \  3.5042   &855\ \ \      1.1386\\
$R_{4,19}$         & 866\ \ \ 0.7500   &{\bf 62\ \ 0.0577}     &626\ \ \ 0.6308     &163\ \ \  1.9576   &632\ \ \      1.0034\\
$R_{6,1}$          & 257\ \ \ 0.2888   &{\bf 65\ \ 0.0603}     &185\ \ \ 0.1803     & 240\ \ \  3.0334   &367\ \ \      0.5613\\
$R_{6,2}$          & 528\ \ \ 0.4939   &{\bf 90\ \ 0.0814}     &543\ \ \ 0.5102     & 285\ \ \  3.5599   &638\ \ \      0.9902\\
$R_{6,3}$          & 172\ \ \ 0.2253   &{\bf 53\ \ 0.0521}     &159\ \ \ 0.2301     & 204\ \ \  2.0102   &299\ \ \      0.4097\\
$R_{6,4}$          & 426\ \ \ 0.4216   &{\bf 90\ \ 0.0792}     &405\ \ \ 0.3529     & 248\ \ \  2.8871   &535\ \ \      0.8490\\
$R_{6,5}$          & 185\ \ \ 0.2308   &{\bf 65\ \ 0.0594}     &116\ \ \ 0.1308     & 226\ \ \  2.3169    &284\ \ \      0.4334\\
\hline
\end{tabular}
\end{table}

\begin{table}[H]
\centering
\caption{The numerical results for multilinear PageRank with $\theta=0.99$.}
\begin{tabular}{cccccc}
\hline
                   &HOPM     &QEHOPM        &RHOPM    &IO      &S \\
                   & IT\ \ CPU  & IT\ \ CPU & IT\ \ CPU  & IT\ \ CPU   & IT\ \ CPU\\
\hline
$R_{3,1}$          & -\ \ \ \ \ \ -          &{\bf 13\ \ 0.0185}       &47\ \ \ 0.0836     & 54\ \ \  0.6360    &52\ \ \      0.0860\\
$R_{3,2}$          & 34 \ \ 0.1222     &{\bf 13\ \  0.0249}      &14\ \ \ 0.0306     & 60\ \ \  0.5520    &56\ \ \      0.0899\\
$R_{3,3}$          & 34\ \ \ 0.1017    &{\bf 13\ \  0.0222}      &14\ \ \ 0.0326     & 60\ \ \  0.5866    &56\ \ \      0.1036\\
$R_{3,4}$          & 148\ \ \ 0.2185   &{\bf 17\ \  0.0248}      &21\ \ \ 0.0369     &32\ \ \   0.4296    &21\ \ \     0.0447\\
$R_{3,5}$          & 27\ \ \ 0.1040    &{\bf 19\ \ 0.0252}      &31\ \ \ 0.0835      &50\ \ \  0.5139     &47\ \ \     0.0887\\
$R_{4,1}$          & -\ \ \ \ \ \ -    &{\bf 70\ \ 0.0438}      &-\ \ \ \ \ \ -      &667\ \ \  6.1107    &-\ \ \ \ \ \ - \\
$R_{4,2}$          & -\ \ \ \ \ \ -    &{\bf 30\ \ 0.0294}      &-\ \ \ \ \ \ -            &283\ \ \  2.0479    &-\ \ \ \ \ \ -  \\
$R_{4,3}$          & -\ \ \ \ \ \ -    &{\bf 32\ \  0.0297}    &-\ \ \ \ \ \ -             &-\ \ \ \ \ \ -      &-\ \ \ \ \ \ - \\
$R_{4,4}$          &-\ \ \ \ \ \ -     &{\bf 29\ \ 0.0299}      &-\ \ \ \ \ \ -            &584\ \    4.5101    &-\ \ \ \ \ \ -\\
$R_{4,5}$          &-\ \ \ \ \ \ -     &{\bf 26\ \  0.0272}     &-\ \ \ \ \ \ -            &386\ \ \  3.5180    &-\ \ \ \ \ \ -\\
$R_{4,6}$          &-\ \ \ \ \ \ -     &{\bf 22\ \  0.0243}     &-\ \ \ \ \ \ -            &276\ \ \  2.5554    &-\ \ \ \ \ \ -\\
$R_{4,7}$          &-\ \ \ \ \ \ -     &{\bf 25\ \  0.0259}      &-\ \ \ \ \ \ -             &427\ \ \  4.9250  &-\ \ \ \ \ \ -\\
$R_{4,8}$          &-\ \ \ \ \ \ -     &{\bf 30\ \  0.0295}     &-\ \ \ \ \ \ -             &191\ \ \  1.8877   &-\ \ \ \ \ \ -\\
$R_{4,9}$          &-\ \ \ \ \ \ -     &{\bf 42\ \  0.0399}      &-\ \ \ \ \ \ -            &391\ \ \  5.5578   &-\ \ \ \ \ \ - \\
$R_{4,10}$         &-\ \ \ \ \ \ -     &{\bf 63\ \ 0.0553}      &-\ \ \ \ \ \ -             &660\ \ \  6.1657   &-\ \ \ \ \ \ - \\
$R_{4,11}$         &-\ \ \ \ \ \ -     &{\bf 62\ \  0.0678}    &-\ \ \ \ \ \ -            &-\ \ \ \ \ \ -       &-\ \ \ \ \ \ -   \\
$R_{4,12}$         &-\ \ \ \ \ \ -     &{\bf 172\ \  0.1404}    &-\ \ \ \ \ \ -             &-\ \ \ \ \ \ -      &-\ \ \ \ \ \ -  \\
$R_{4,13}$         &-\ \ \ \ \ \ -     &{\bf 113\ \ 0.1121}      &-\ \ \ \ \ \ -            & 623\ \ \  5.8067   &-\ \ \ \ \ \ -   \\
$R_{4,14}$         &-\ \ \ \ \ \ -     &{\bf 34\ \ 0.0348}      &-\ \ \ \ \ \ -             & 334\ \ \  3.3455   &-\ \ \ \ \ \ -   \\
$R_{4,15}$         &-\ \ \ \ \ \ -     &{\bf 34\ \  0.0332}     &-\ \ \ \ \ \ -             & 283\ \ \  2.9808   &-\ \ \ \ \ \ -    \\
$R_{4,16}$         &-\ \ \ \ \ \ -     &{\bf 29\ \  0.0310}     &-\ \ \ \ \ \ -             & 638\ \ \  6.8728   &-\ \ \ \ \ \ - \\
$R_{4,17}$         &-\ \ \ \ \ \ -     &{\bf 121\ \ 0.1012}     &510\ \ \ 0.1209      &-\ \ \ \ \ \ -      &642\ \ \ 0.9793\\
$R_{4,18}$         &-\ \ \ \ \ \ -     &{\bf 229\ \ 0.1963}     &-\ \ \ \ \ \ -             & 529\ \ \  5.5947   &-\ \ \ \ \ \ - \\
$R_{4,19}$         &-\ \ \ \ \ \ -     &{\bf 42\ \  0.0390}     &-\ \ \ \ \ \ -             & 185\ \ \  1.3840   &-\ \ \ \ \ \ -\\
$R_{6,1}$          &-\ \ \ \ \ \ -     &{\bf 41\ \  0.0310}     &-\ \ \ \ \ \ -             & 902\ \ \  9.9988   &-\ \ \ \ \ \ - \\
$R_{6,2}$          &-\ \ \ \ \ \ -     &{\bf 76\ \  0.0599}     &-\ \ \ \ \ \ -             & 765\ \ \  8.1938   &-\ \ \ \ \ \ - \\
$R_{6,3}$          &-\ \ \ \ \ \ -     &{\bf 320\ \ 0.2306}     &-\ \ \ \ \ \ -        &-\ \ \ \ \ \ -           &-\ \ \ \ \ \ -  \\
$R_{6,4}$          &-\ \ \ \ \ \ -     &{\bf 74\ \  0.0557}     &-\ \ \ \ \ \ -             & 420\ \ \  5.3745   &-\ \ \ \ \ \ -   \\
$R_{6,5}$          &-\ \ \ \ \ \ -     &{\bf 404\ \  0.3316}    &-\ \ \ \ \ \ -             &-\ \ \ \ \ \ -            &-\ \ \ \ \ \ -   \\
\hline
\end{tabular}
\end{table}

As we can see from Tables 2-6, the proposed QEHOPM method is faster than HOPM, RHOPM, IO and S methods. The QEHOPM method spends less iterations or CPU time than that of HOPM, RHOPM, IO and S methods. In particular, QEHOPM method has the most reliable convergence, which can solve successfully all of the benchmark set of 29 problems even when $\theta=0.99$. In a word, from the above results, we conclude that QEHOPM algorithm is very effective and competitive.


\section{Conclusion}
In this paper, we have proposed three accelerated higher-order power method for higher-order Markov chains and multilinear PageRank, referred to as the HOPMM-I/HOPMM-II and QEHOPM, respectively. In particular, the QEHOPM method is non-parametric. We established the convergence results for the proposed algorithms. Numerical experiments are carried out to illustrate that HOPMM-I/HOPMM-II outperform the higher-order power method and the QEHOPM is the best one.

\section*{Appendix}
\noindent\textbf{Test Examples}: The first three tensors come from DNA sequence data in the works of Raftery et al. \cite{Raftery94}. On the other hand, their orders $m$ are 3 and their numbers of states $n$ are 3 or 4 by considering three categories(\{A/G,C,T). By using the Matlab multi-dimensional array notation, the transition probability tensors are given as follows.\\

(i)$\mathcal{P}(:,:,1)=
\left(\begin{array}{ccc}
0.6000 & 0.4083 & 0.4935\\
0.2000 &0.2568 & 0.2426\\
0.2000 & 0.3349  &  0.2639\\
\end{array}\right),$
 $\mathcal{P}(:,:,2)=
\left(\begin{array}{cccc}
0.5217 & 0.3300 & 0.4152\\
0.2232 &0.2800 & 0.2658\\
0.2551& 0.3900  &  0.3190\\
\end{array}\right),$

 $$\mathcal{P}(:,:,3)=
\left(\begin{array}{cccc}
0.5565 & 0.3648 & 0.4500\\
0.2174&0.2742 & 0.2600\\
0.2261 & 0.3610  &  0.2900\\
\end{array}\right).$$
(ii)$\mathcal{P}(:,:,1)=
\left(\begin{array}{ccc}
0.5200 & 0.2986 & 0.4462\\
0.2700 &0.3930 & 0.3192\\
0.2100 & 0.3084  &  0.2346\\
\end{array}\right),$
 $\mathcal{P}(:,:,2)=
\left(\begin{array}{cccc}
0.6514 & 0.4300 & 0.5776\\
0.1970 &0.3200 & 0.2462\\
0.1516& 0.2500  &  0.1762\\
\end{array}\right),$

 $$\mathcal{P}(:,:,3)=
\left(\begin{array}{cccc}
0.5638 & 0.3424 & 0.4900\\
0.2408&0.3638 & 0.2900\\
0.1954 & 0.2938  &  0.2200\\
\end{array}\right).$$

(iii)$$\mathcal{P}(:,:,1)=
\left(\begin{array}{cccc}
0.2091 & 0.2834 & 0.2194 & 0.1830\\
0.3371 &0.3997 & 0.3219& 0.3377\\
0.3265 & 0.0560  &  0.3119& 0.2961\\
0.1723 & 0.2608  &  0.1468& 0.1832\\
\end{array}\right),$$

$$\mathcal{P}(:,:,2)=
\left(\begin{array}{cccc}
0.1952 & 0.2695 & 0.2055& 0.1690\\
0.3336 &0.3962 & 0.3184& 0.3342\\
0.2954& 0.0249  &  0.2808& 0.2650\\
0.1758& 0.3094  &  0.1953& 0.2318\\
\end{array}\right),$$

 $$\mathcal{P}(:,:,3)=
\left(\begin{array}{cccc}
0.3145 & 0.3887 & 0.3248& 0.2883\\
0.0603&0.1203 & 0.0451& 0.0609\\
0.2293 & 0.3628  &  0.2487& 0.2852\\
0.2293 & 0.3628  &  0.2487& 0.2852\\
\end{array}\right).$$

 $$\mathcal{P}(:,:,4)=
\left(\begin{array}{cccc}
0.1685 & 0.2429 & 0.1789& 0.1425\\
0.3553&0.4180 & 0.3402& 0.3559\\
0.3189 & 0.0484  &  0.3043& 0.2885\\
0.1571 & 0.2907  &  0.1766& 0.2131\\
\end{array}\right).$$

By considering three categories (\{A,C/T,G\}), we construct a transition probability tensor of order 4 and dimension 3 for the DNA sequence in \cite{NCBI}:\\
(iv)$\mathcal{P}(:,:,1,1)=
\left(\begin{array}{ccc}
0.3721 & 0.2600 & 0.4157\\
0.4477 &0.5000 & 0.4270\\
0.1802 & 0.2400  &  0.1573\\
\end{array}\right),$
$\mathcal{P}(:,:,2,1)=
\left(\begin{array}{cccc}
0.3692 & 0.2673 & 0.3175\\
0.4667 &0.5594 & 0.5079\\
0.1641& 0.1733  &  0.1746\\
\end{array}\right),$

 $\mathcal{P}(:,:,3,1)=
\left(\begin{array}{cccc}
0.4227 & 0.2958 & 0.2353\\
0.4124&0.5563 & 0.5588\\
0.1649 & 0.1479  &  0.2059\\
\end{array}\right)$,
$\mathcal{P}(:,:,1,2)=
\left(\begin{array}{cccc}
0.3178 & 0.2632 & 0.3194\\
0.5212&0.6228 & 0.5833\\
0.1610 & 0.1140  &  0.0972\\
\end{array}\right)$,

 $\mathcal{P}(:,:,2,2)=
\left(\begin{array}{cccc}
0.2836 & 0.2636 & 0.3042\\
0.5012&0.6000 & 0.5250\\
0.2152 & 0.1364  &  0.1708\\
\end{array}\right)$,
$\mathcal{P}(:,:,3,2)=
\left(\begin{array}{cccc}
0.3382 & 0.2396 & 0.3766\\
0.5147&0.6406 & 0.4935\\
0.1471 & 0.1198  &  0.1299\\
\end{array}\right)$,

 $\mathcal{P}(:,:,1,3)=
\left(\begin{array}{cccc}
0.3204 & 0.2985 & 0.3500\\
0.4854&0.5000 & 0.5000\\
0.1942 & 0.2015  &  0.1500\\
\end{array}\right)$,
 $\mathcal{P}(:,:,2,3)=
\left(\begin{array}{cccc}
0.4068 & 0.2816 & 0.3594\\
0.3898&0.5143& 0.4219\\
0.2034 & 0.2041  &  0.2188\\
\end{array}\right)$,

 $$\mathcal{P}(:,:,3,3)=
\left(\begin{array}{cccc}
0.3721 & 0.3529 & 0.3000\\
0.5349&0.3971 & 0.5500\\
0.0930 & 0.2500  &  0.1500\\
\end{array}\right).$$

\section*{Acknowledgement}
This work was supported in part by the National Natural Science Foundation of China (No.11661007, 11701097, and 61806004), in part by Natural Science Foundation of Zhejiang Province (No.LD19A010002), and in part by Natural Science Foundation of Jiangxi Province (20181BAB201007).


\end{document}